\documentclass[prodmode,acmjacm,oneside]{acmsmall}

\usepackage{bm}
\usepackage{latexsym,xspace,calc,amsthm}
\usepackage{amsmath,amssymb} 

\makeatletter
\newcommand{\domsub}{\mathrel{\mathpalette\dom@sub\relax}}
\newcommand{\dom@sub}[2]{%
  \ooalign{%
    $\m@th#1\triangleright$\cr
    \hidewidth$\m@th#1\text{\scalebox{.65}[1]{$-$}}$\hidewidth\cr
  }%
}
\newcommand{\ransub}{\mathrel{\mathpalette\ran@sub\relax}}
\newcommand{\ran@sub}[2]{%
  \ooalign{%
    $\m@th#1\triangleleft$\cr
    \hidewidth$\m@th#1\text{\scalebox{.65}[1]{$-$}}$\hidewidth\cr
  }%
}
\makeatother

\makeatletter
\def\@biblabel#1{[#1]} 
\def\thebibliography#1{%
    \footnotesize
    \refsection*{{\refname}
        \@mkboth{\uppercase{\refname}}{\uppercase{\refname}}%
    }
    \list{\@biblabel{\@arabic\c@enumiv}}
       {\settowidth\labelwidth{\@biblabel{#1}}%
        \leftmargin\labelwidth
        \advance\leftmargin\bibindent
        \itemindent-\bibindent
        \itemsep2pt
        \parsep \z@
        \usecounter{enumiv}
        \let\p@enumiv\@empty
        \renewcommand\theenumiv{\@arabic\c@enumiv}%
    }%
    \let\newblock\@empty
    \sloppy
    \sfcode`\.=1000\relax
}
\makeatother

\usepackage{enumerate,mdwlist} 

\makecompactlist{enumeratei*}{enumeratei}
\usepackage{makeidx}

\usepackage{xparse}
\usepackage{tikz}
\usepackage{tikz-cd}
\usetikzlibrary{quotes}

\usetikzlibrary{calc}
\DeclareDocumentCommand{\hcancel}{mO{0pt}O{1pt}O{0pt}O{-1pt}}{%
    \tikz[baseline=(tocancel.base)]{
        \node[inner sep=0pt,outer sep=0pt] (tocancel) {#1};
        \draw[gray] ($(tocancel.south west)+(#2,#3)$) -- ($(tocancel.north east)+(#4,#5)$);
    }%
}%

\usepackage{microtype}
\usepackage{tgtermes}
\usepackage{fix-cm}

\usepackage{amssymb,stackengine,graphicx,trimclip}
\stackMath
\newcommand\frightarrow{\scalebox{.4}[.3]{%
  \clipbox{0pt 0pt 2pt 0pt}{$\blacktriangleright$}\kern1.8pt}}
\newcommand\fleftarrow{\scalebox{.4}[.3]{%
  \kern1.8pt\clipbox{2pt 0pt 0pt 0pt}{$\blacktriangleleft$}}}
\newcommand\subsetrleads{\mathrel{%
  \stackengine{-1.1pt}{\subseteq}{\frightarrow\kern.3pt}{U}{r}{F}{T}{S}}}

\newcommand\subsetlleads{\mathrel{%
  \stackengine{-1.1pt}{\subseteq}{\kern.4pt\fleftarrow}{U}{l}{F}{T}{S}}}

\newcommand\subsetsim{\mathrel{%
  \ooalign{\raise0.2ex\hbox{$\subset$}\cr\hidewidth\raise-0.8ex\hbox{\scalebox{0.9}{$\sim$}}\hidewidth\cr}}}

\usepackage{float}
\floatstyle{boxed} \restylefloat{figure}

\usepackage[all]{xy}

\DeclareSymbolFont{sfletters}{OML}{cmbrm}{m}{it}
\DeclareMathSymbol{\salpha}{\mathord}{sfletters}{"0B}
\DeclareMathSymbol{\sbeta}{\mathord}{sfletters}{"0C}
\DeclareMathSymbol{\sgamma}{\mathord}{sfletters}{"0D}
\DeclareMathSymbol{\sdelta}{\mathord}{sfletters}{"0E}
\DeclareMathSymbol{\sepsilon}{\mathord}{sfletters}{"0F}
\DeclareMathSymbol{\szeta}{\mathord}{sfletters}{"10}
\DeclareMathSymbol{\seta}{\mathord}{sfletters}{"11}
\DeclareMathSymbol{\stheta}{\mathord}{sfletters}{"12}
\DeclareMathSymbol{\siota}{\mathord}{sfletters}{"13}
\DeclareMathSymbol{\skappa}{\mathord}{sfletters}{"14}
\DeclareMathSymbol{\slambda}{\mathord}{sfletters}{"15}
\DeclareMathSymbol{\smu}{\mathord}{sfletters}{"16}
\DeclareMathSymbol{\snu}{\mathord}{sfletters}{"17}
\DeclareMathSymbol{\sxi}{\mathord}{sfletters}{"18}
\DeclareMathSymbol{\spi}{\mathord}{sfletters}{"19}
\DeclareMathSymbol{\srho}{\mathord}{sfletters}{"1A}
\DeclareMathSymbol{\ssigma}{\mathord}{sfletters}{"1B}
\DeclareMathSymbol{\stau}{\mathord}{sfletters}{"1C}
\DeclareMathSymbol{\supsilon}{\mathord}{sfletters}{"1D}
\DeclareMathSymbol{\sphi}{\mathord}{sfletters}{"1E}
\DeclareMathSymbol{\schi}{\mathord}{sfletters}{"1F}
\DeclareMathSymbol{\spsi}{\mathord}{sfletters}{"20}
\DeclareMathSymbol{\somega}{\mathord}{sfletters}{"21}
\DeclareMathSymbol{\svarepsilon}{\mathord}{sfletters}{"22}

\DeclareRobustCommand{\lowercasei}{i}

\makeatletter
\def\pmb@#1#2{\setbox8\hbox{$\m@th#1{#2}$}%
  \setboxz@h{$\m@th#1\mkern-.1mu$}\pmbraise@\wdz@
  \binrel@{#2}%
  \dimen@-\wd8 %
  \binrel@@{%
    \mkern-.1mu\copy8 %
    \kern\dimen@\mkern-.2mu\copy8 %
    \kern\dimen@\mkern-.3mu\copy8 %
    \kern\dimen@\mkern-.4mu\copy8 %
    \kern\dimen@\mkern.1mu\copy8 %
    \kern\dimen@\mkern.2mu\copy8 %
    \kern\dimen@\mkern.3mu\copy8 %
    \kern\dimen@\mkern.0mu\raise\pmbraise@\copy8 %
    \kern\dimen@\mkern.4mu\box8 %
           }%
}
\makeatother

\newcommand\Cmp[1]{\rulefont{Cmp_{\mathit{#1}}}}

\newcommand\nfin{\mathbin{\varepsilon}}

\newcommand\fv{\f{fv}}
\newcommand\consts{\f{const}}

\newcommand\upsome[1]{{}^{{\plus}#1}}
\newcommand\downsome[1]{{}^{{\downarrow}}}

\newcommand\downone[1]{{\downarrow}^{\hspace{-1.4pt}#1}}

\newcommand\card{\myfresh}

\let\myfresh\#
\def\#{\ensuremath{\text{\tt\myfresh}}}
\newcommand\supports[0]{\mathop{\texttt{\it\$}}} 

\newcommand\univ{\tf W}
\newcommand\univi[1]{\univ_{#1}}
\newcommand\hfs[1]{\tf{HFS}_{#1}}

\newcommand\finsubseteq{\mathbin{\subseteq_{\text{\it fin}}}}
\newcommand\level{\f{lev}}
\newcommand\compressthis[1]{\pmb{\hspace{.8pt}\raisebox{.5pt}{\scalebox{.85}{$#1$}}\hspace{.2pt}}}

\newcommand\tneg{{\pmb\neg}}

\newcommand\tbot{{\compressthis\bot}}
\newcommand\teq{{\pmb{\text{=}}}}
\newcommand\tand{{\compressthis\wedge}}
\newcommand\tor{{\compressthis\vee}}
\newcommand\timp{{\compressthis\Rightarrow}}

\newcommand\tiff{\compressthis{\Leftrightarrow}} 

\newcommand\tin{{\compressthis{\hspace{0.3pt}\varepsilon\hspace{0.5pt}}}} 
\newcommand\tst[2]{\compressthis{\{}#1\hspace{.25pt}\compressthis{\mid}\hspace{.4pt}#2\compressthis{\}}}
\newcommand\stxst[2]{{\tst{#1}{#2}}}

\newcommand\tall{{\compressthis{\forall}}}
\newcommand\texi{{\compressthis{\exists}}}
\newcommand\compactwedge[1]{\hspace{-1ex}\bigwedge_{#1}\hspace{-1em}}

\renewcommand\land{\wedge}

\newcommand\limp{\Rightarrow}

\newcommand\minus{{\text{-}}}
\newcommand\plus{{\text{+}}}

\usepackage{upgreek}

 \renewenvironment{thebibliography}[1]{%
   \begin{odlthebibliography}{#1}%
     \setlength{\parskip}{0ex}%
     \setlength{\itemsep}{3pt}%
     \fontsize{10}{10} 
     \selectfont
}%
 {%
   \end{odlthebibliography}%
 }

\usepackage[colorlinks]{hyperref}
\hypersetup{urlcolor=black,colorlinks}

\usepackage{fancybox}
\newlength{\mylength}
{\setlength{\fboxsep}{5pt}
\setlength{\mylength}{\linewidth}%
\addtolength{\mylength}{-2\fboxsep}%
\addtolength{\mylength}{-2\fboxrule}%
\Sbox
\minipage{\mylength}%
\setlength{\abovedisplayskip}{0pt}%
\setlength{\belowdisplayskip}{0pt}%
$$}%
{$$\endminipage\endSbox
{\setlength{\abovedisplayskip}{1pt}%
\setlength{\belowdisplayskip}{0pt}%
\[\fbox{\TheSbox}\]}}

{\setlength{\fboxsep}{5pt}
\setlength{\mylength}{\linewidth}%
\addtolength{\mylength}{-2\fboxsep}%
\addtolength{\mylength}{-2\fboxrule}%
\Sbox
\minipage{\mylength}%
\setlength{\abovedisplayskip}{5pt}%
\setlength{\belowdisplayskip}{5pt}%
}%
{\endminipage\endSbox
{\setlength{\abovedisplayskip}{1pt}%
\setlength{\belowdisplayskip}{0pt}%
\[\fbox{\TheSbox}\]}}

\usepackage{graphics}
\usepackage{ifthen}
\usepackage{verbatim}
\newdimen\proofrulebreadth \proofrulebreadth=.05em
\newdimen\proofdotseparation \proofdotseparation=1.25ex
\newdimen\proofrulebaseline \proofrulebaseline=2ex
\newcount\proofdotnumber \proofdotnumber=3
\let\then\relax
\def\hfi{\hskip0pt plus.0001fil}
\mathchardef\squigto="3A3B
\newif\ifinsideprooftree\insideprooftreefalse
\newif\ifonleftofproofrule\onleftofproofrulefalse
\newif\ifproofdots\proofdotsfalse
\newif\ifdoubleproof\doubleprooffalse
\let\wereinproofbit\relax
\newdimen\shortenproofleft
\newdimen\shortenproofright
\newdimen\proofbelowshift
\newbox\proofabove
\newbox\proofbelow
\newbox\proofrulename
\def\shiftproofbelow{\let\next\relax\afterassignment\setshiftproofbelow\dimen0 }
\def\shiftproofbelowneg{\def\next{\multiply\dimen0 by-1 }%
\afterassignment\setshiftproofbelow\dimen0 }
\def\setshiftproofbelow{\next\proofbelowshift=\dimen0 }
\def\setproofrulebreadth{\proofrulebreadth}

\def\prooftree{
\ifnum  \lastpenalty=1
\then   \unpenalty
\else   \onleftofproofrulefalse
\fi
\ifonleftofproofrule
\else   \ifinsideprooftree
        \then   \hskip.5em plus1fil
        \fi
\fi
\bgroup
\setbox\proofbelow=\hbox{}\setbox\proofrulename=\hbox{}%
\let\justifies\proofover\let\leadsto\proofoverdots\let\Justifies\proofoverdbl
\let\using\proofusing\let\[\prooftree
\ifinsideprooftree\let\]\endprooftree\fi
\proofdotsfalse\doubleprooffalse
\let\thickness\setproofrulebreadth
\let\shiftright\shiftproofbelow \let\shift\shiftproofbelow
\let\shiftleft\shiftproofbelowneg
\let\ifwasinsideprooftree\ifinsideprooftree
\insideprooftreetrue
\setbox\proofabove=\hbox\bgroup$\displaystyle 
\let\wereinproofbit\prooftree
\shortenproofleft=0pt \shortenproofright=0pt \proofbelowshift=0pt
\onleftofproofruletrue\penalty1
}

\def\eproofbit{
\ifx    \wereinproofbit\prooftree
\then   \ifcase \lastpenalty
        \then   \shortenproofright=0pt  
        \or     \unpenalty\hfil         
        \or     \unpenalty\unskip       
        \else   \shortenproofright=0pt  
        \fi
\fi
\global\dimen0=\shortenproofleft
\global\dimen1=\shortenproofright
\global\dimen2=\proofrulebreadth
\global\dimen3=\proofbelowshift
\global\dimen4=\proofdotseparation
\global\count255=\proofdotnumber
$\egroup  
\shortenproofleft=\dimen0
\shortenproofright=\dimen1
\proofrulebreadth=\dimen2
\proofbelowshift=\dimen3
\proofdotseparation=\dimen4
\proofdotnumber=\count255
}

\def\proofover{
\eproofbit 
\setbox\proofbelow=\hbox\bgroup 
\let\wereinproofbit\proofover
$\displaystyle
}%
\def\proofoverdbl{
\eproofbit 
\doubleprooftrue
\setbox\proofbelow=\hbox\bgroup 
\let\wereinproofbit\proofoverdbl
$\displaystyle
}%
\def\proofoverdots{
\eproofbit 
\proofdotstrue
\setbox\proofbelow=\hbox\bgroup 
\let\wereinproofbit\proofoverdots
$\displaystyle
}%
\def\proofusing{
\eproofbit 
\setbox\proofrulename=\hbox\bgroup 
\let\wereinproofbit\proofusing
\kern0.3em$
}

\def\endprooftree{
\eproofbit 
  \dimen5 =0pt
\dimen0=\wd\proofabove \advance\dimen0-\shortenproofleft
\advance\dimen0-\shortenproofright
\dimen1=.5\dimen0 \advance\dimen1-.5\wd\proofbelow
\dimen4=\dimen1
\advance\dimen1\proofbelowshift \advance\dimen4-\proofbelowshift
\ifdim  \dimen1<0pt
\then   \advance\shortenproofleft\dimen1
        \advance\dimen0-\dimen1
        \dimen1=0pt
        \ifdim  \shortenproofleft<0pt
        \then   \setbox\proofabove=\hbox{%
                        \kern-\shortenproofleft\unhbox\proofabove}%
                \shortenproofleft=0pt
        \fi
\fi
\ifdim  \dimen4<0pt
\then   \advance\shortenproofright\dimen4
        \advance\dimen0-\dimen4
        \dimen4=0pt
\fi
\ifdim  \shortenproofright<\wd\proofrulename
\then   \shortenproofright=\wd\proofrulename
\fi
\dimen2=\shortenproofleft \advance\dimen2 by\dimen1
\dimen3=\shortenproofright\advance\dimen3 by\dimen4
\ifproofdots
\then
        \dimen6=\shortenproofleft \advance\dimen6 .5\dimen0
        \setbox1=\vbox to\proofdotseparation{\vss\hbox{$\cdot$}\vss}%
        \setbox0=\hbox{%
                \advance\dimen6-.5\wd1
                \kern\dimen6
                $\vcenter to\proofdotnumber\proofdotseparation
                        {\leaders\box1\vfill}$%
                \unhbox\proofrulename}%
\else   \dimen6=\fontdimen22\the\textfont2 
        \dimen7=\dimen6
        \advance\dimen6by.5\proofrulebreadth
        \advance\dimen7by-.5\proofrulebreadth
        \setbox0=\hbox{%
                \kern\shortenproofleft
                \ifdoubleproof
                \then   \hbox to\dimen0{%
                        $\mathsurround0pt\mathord=\mkern-6mu%
                        \cleaders\hbox{$\mkern-2mu=\mkern-2mu$}\hfill
                        \mkern-6mu\mathord=$}%
                \else   \vrule height\dimen6 depth-\dimen7 width\dimen0
                \fi
                \unhbox\proofrulename}%
        \ht0=\dimen6 \dp0=-\dimen7
\fi
\let\doll\relax
\ifwasinsideprooftree
\then   \let\VBOX\vbox
\else   \ifmmode\else$\let\doll=$\fi
        \let\VBOX\vcenter
\fi
\VBOX   {\baselineskip\proofrulebaseline \lineskip.2ex
        \expandafter\lineskiplimit\ifproofdots0ex\else-0.6ex\fi
        \hbox   spread\dimen5   {\hfi\unhbox\proofabove\hfi}%
        \hbox{\box0}%
        \hbox   {\kern\dimen2 \box\proofbelow}}\doll%
\global\dimen2=\dimen2
\global\dimen3=\dimen3
\egroup 
\ifonleftofproofrule
\then   \shortenproofleft=\dimen2
\fi
\shortenproofright=\dimen3
\onleftofproofrulefalse
\ifinsideprooftree
\then   \hskip.5em plus 1fil \penalty2
\fi
}

\newcommand\ns{\mathsf}

\def\:{{\hspace{-1pt}{:}\hspace{-1.25pt}{:}\hspace{-.5pt}}}
\usepackage{relsize}

\newcommand\id{\f{id}}

\newcommand\ssm{{{:}\text{=}}}

\newcommand\deffont[1]{{\bfseries #1}}

\newcommand\minussign{{{\text{-}}}}
\newcommand\mone{{{\minussign 1}}}
\newcommand\liff{\Leftrightarrow}

\newcommand\supp{\f{supp}}
\newcommand\adjsupp{\f{supp}_{\hspace{-.7pt}\pmb{\mact}}}

\newcommand\f[1]{\mathit{#1}}
\newcommand\tf[1]{\mathsf{#1}}

\newcommand\atoms{\ensuremath{\mathbb{A}}\xspace}
\newcommand\btoms{\ensuremath{\mathbb{B}}\xspace}

\newcommand\varsi[1]{\tf{Var}_{\hspace{-.7pt}{#1}}}

\newcommand\powerset{\f{pow}}
\newcommand\powersetfs{\powerset_{\hspace{-1pt}\f{fs}}}

\newcommand\lmodel{[\hspace{-0.2em}[}
\newcommand\rmodel{]\hspace{-0.2em}]}
\newcommand\model[1]{{\lmodel #1 \rmodel}}

\newcommand\act[0]{{\cdot}}
\newcommand\mact[0]{{{\cdot}{\cdot}}}

\newcommand{\defeq}
  {\stackrel{\mathrm{def}}{\,=\,}}

\newcommand\fix{\f{fix}}

\newcommand\ment[0]{\mathrel{\vDash}}

\newcounter{jamieitemcounter}

\newtheoremstyle{jamiestyle}
  {4pt}
  {0pt}
  {\it}
  {0pt}
  {\sc}
  {.}
  { }
  {}
\theoremstyle{jamiestyle}
\newtheorem{thrm}{Theorem}[subsection]
\newtheorem{prop}[thrm]{Proposition}
\newtheorem{lemm}[thrm]{Lemma}
\newtheorem{corr}[thrm]{Corollary}
\newtheoremstyle{jamienfstyle}
  {4pt}
  {0pt}
  {\normalfont}
  {0pt}
  {\sc}
  {.}
  { }
  {}
\theoremstyle{jamienfstyle}
\newtheorem{nttn}[thrm]{Notation}
\newtheorem{defn}[thrm]{Definition}
\newtheorem{xmpl}[thrm]{Example}
\newtheorem{rmrk}[thrm]{Remark}

\usepackage{stmaryrd}

\newcommand\rulefont[1]{\ensuremath{{\mathrm{\bf (#1)}}}}

\newcommand\Forall[1]{\forall #1.}
\newcommand\Exists[1]{\exists #1.}

\usepackage{datetime}

\makeatletter
\@ifpackageloaded{fdsymbol}\@tempswafalse\@tempswatrue
\if@tempswa
  \newcommand{\fdsy@scale}{1.0}
  \newcommand\fdsy@mweight@normal{Book}
  \newcommand\fdsy@mweight@small{Book}
  \newcommand\fdsy@bweight@normal{Medium}
  \newcommand\fdsy@bweight@small{Medium}
  \DeclareFontFamily{U}{FdSymbolA}{}
  \DeclareSymbolFont{fdsymbols}{U}{FdSymbolA}{m}{n}%
  \SetSymbolFont{symbols}{bold}{U}{FdSymbolA}{b}{n}%
  \DeclareFontShape{U}{FdSymbolA}{m}{n}{
      <-7.1> s * [\fdsy@scale] FdSymbolA-\fdsy@mweight@small
      <7.1-> s * [\fdsy@scale] FdSymbolA-\fdsy@mweight@normal
  }{}
  \DeclareFontShape{U}{FdSymbolA}{b}{n}{
      <-7.1> s * [\fdsy@scale] FdSymbolA-\fdsy@bweight@small
      <7.1-> s * [\fdsy@scale] FdSymbolA-\fdsy@bweight@normal
  }{}
  \DeclareMathSymbol{\aleph}{\mathord}{fdsymbols}{"C7}
  \DeclareMathSymbol{\beth}{\mathord}{fdsymbols}{"C8}
  \DeclareMathSymbol{\gimel}{\mathord}{fdsymbols}{"C9}
  \DeclareMathSymbol{\daleth}{\mathord}{fdsymbols}{"CA}
\fi
\makeatother

\title{Consistency of Quine's New Foundations}
\author{\href{http://www.gabbay.org.uk}{Murdoch J. Gabbay}
\affil{\href{http://www.gabbay.org.uk}{\it http://www.gabbay.org.uk}}
}

\begin{abstract}
We provide an elementary consistency proof of Quine's New Foundations, by a construction using interated nominal powersets.
\end{abstract}

\makeindex

\begin{document}
\maketitle
\tableofcontents

\section{Introduction}
\label{sect.intro}

\subsection{New foundations}

Consider the following \emph{invalid} reasoning: define $x=\{a\mid a{\not\in} a\}$.
It is easy to check that $x\in x$ if and only if $x\not\in x$.
This is Russell's paradox and is one of the central paradoxes of (naive) set theory.

Zermelo-Fraenkel set theory (\deffont{ZF}) avoids paradox by insisting instead that $a$ be \emph{guarded}; we can only form $\{a{\in}y\mid a{\not\in} a\}$ where $y$ is already known to be a set.
The price we pay for this is that we cannot form `reasonable' sets such as the \deffont{universal set} $\{a\mid \top\}$ (the set of all objects) or the set of `all sets with 2 elements', and so on.
In ZF, these are \emph{proper classes}.\footnote{A nice historical account of Russell's paradox is in \cite{griffin:prehrp}. For ZF set theory, see e.g. \cite{jech:sett}.}

Quine's New Foundations (\deffont{NF}) avoids paradox by restricting to stratifiable comprehensions~\cite{quine:newfml},
such that every variable and term can be assigned a \emph{level} and $t\in s$ only appears in the $\phi$ of a comprehension $\{a\mid \phi\}$ when $\level(s)=\level(t)\plus 1$. 
For example, the comprehensions $\{a\mid a\in a\}$ and $\{a\mid a\not\in a\}$ are not stratifiable, because whatever level $i$ we assign to $a$, we cannot make $i$ be equal to $i\plus 1$.

We can stratify $\{a\mid\top\}$ so we can still form the universal set in NF (and `has 2 elements' is also stratifiable).
Excellent discussions are in \cite{forster:settus} and \cite{holmes:elestu}, and a clear summary with a brief but well-chosen bibliography is in \cite{forster:quinfs}.

At the time of writing there is no canonical proof of consistency for NF.\footnote{Holmes has a claimed proof.} 
This has been the situation since NF was introduced in 1937 in~\cite{quine:newfml} --- though a related theory NFU (NF with urelemente) is known consistent by a proof of Jensen~\cite{jensen:oncs,holmes:elestu}.

\begin{rmrk}[Overview of the consistency proof]
\leavevmode
\begin{enumerate}
\item
Section~\ref{sect.intro} is the Introduction.  You Are Here.
\item
In Section~\ref{sect.syntax.and.entailment} we introduce stratified syntax (the syntax of \emph{Typed Set Theory}), define TST+ as a first-order theory, and discuss how this is equiconsistent with NF. 
\item
In Section~\ref{sect.hfs} we build a universe $\hfs{}$ of hereditarily finitely supported sets (\`a la Fraenkel-Mostowski independence proof, or \`a la nominal techniques for syntax-with-binding).
This is essentially off the shelf: the reader can skip it in the first instance.
Note that this paper is self-contained: it should not be necessary to look up prior material on Fraenkel-Mostowksi set theory or nominal techniques, to follow the material below. 
\item
In Section~\ref{sect.map.between} we exploit the finite support condition on $\hfs{}$ to note the existence of bijections $f_i:\hfs{i}\to\hfs{i\plus 1}$.
The proof is by counting cardinalities: in Remark~\ref{rmrk.sketch.proof} we discuss intuitively why it works.
The reader can read this if they like counting arguments, and skip it if they don't.
\item
In Section~\ref{sect.the.model} we build a stratified sets model $\univi{}$. 
It has a particularity that it undergoes a `phase change' at some fixed but arbitrary level $n$: $\univi{i}=\hfs{i}$ for $0\leq i\leq n$, whereas $\univi{n\plus i}$ is in a fairly simple way just a copy of $\univi{n}$, for $i\geq 0$.
Aside from $f$, the meat of the proof is here.
\end{enumerate}
There is plenty of exposition throughout this paper, including
a technical but high-level discussion of the construction in Subsection~\ref{subsect.some.intuitions}, and 
many examples, including Examples~\ref{xmpl.hfs.examples}, \ref{xmpl.updown.examples}, \ref{xmpl.colour}, and a demonstration of the model in action in Example~\ref{xmpl.model}.

The reader who wants to jump to a fulcrum of the consistency construction could start at Definition~\ref{defn.model.tin} (where the $\in$-structure of the model is determined) and work forward through the technical results in Subsection~\ref{subsect.comprehension} (where it is shown that the model satisfies typical ambiguity and extensionality and has enough detail to interpret comprehension), while referring back to the nominal and other dependencies as they are used.
\end{rmrk}

\section{The syntax and entailment of TST+}
\label{sect.syntax.and.entailment}

\subsection{Syntax and entailment}

\begin{rmrk}
\label{rmrk.TZT}
We will work mostly with TST+ rather than NF, for technical reasons.
Intuitively, \emph{Typed Set Theory} TST+ replaces the stratifia\emph{bility} condition of NF, with a \emph{stratification} condition on syntax --- so levels are fixed --- along with a \emph{typical ambiguity} axiom which relates truth at one level with truth at another. 
TST+ and NF are known equiconsistent \cite{specker:typa} (we reproduce the core of the argument in Appendix~\ref{sect.tst.nf} for reference).
\end{rmrk}

\begin{defn}
\label{defn.syntax}
\leavevmode\begin{enumerate}
\item
For each $i\in\mathbb N$ fix a countably infinite set $\varsi{i}$ of \deffont{variable symbols} of level $i$, and a set $\univi{i}$ of \deffont{constants} of level $i$.

If $a\in\varsi{i}$ then define $\level(a)=i$, and if $x\in\univi{i}$ then define $\level(x)=i$.
\item
Define \deffont{raw predicates} and \deffont{raw terms} 
by the following BNF grammar:\footnote{`$\phi\timp\phi$' in this grammar is not a typo: as standard in BNF grammar notation, instances of $\phi$ represent any (not necessarily equal) predicate while also indicating the typical root name(s) of variables typically ranging over the datatype. The $\phi$ and \emph{also} $\psi$ on the left indicates that $\phi$ and \emph{also} $\psi$ will range over predicates.} 
$$
\begin{array}{r@{\ }l}
\phi,\psi ::=& \tbot \mid \phi\timp\phi \mid \tall a.\phi \mid t\tin t 
\\
s,t ::=& a\in\varsi{i} \mid x\in\univi{i} 
\end{array}
$$
\item\label{item.TZT+.syntax}
Define \deffont{stratified predicates} to be those predicates such that 
$t\tin s$ only appears when $\level(t)\plus 1=\level(s)$.
\end{enumerate}
\end{defn}

\begin{rmrk}
The intuition of $y\tin x$ is `$y$ is an element of $x$'.
How this intuition is interpreted in our model, is subtle: see Definition~\ref{defn.model.tin}.

The intuitions of $\tbot$, $\timp$, and $\tall$ are more straightforward: false $\bot$, classical implication $\limp$, and (classical) universal quantification $\forall$ respectively; see Definition~\ref{defn.denotation}.
Our typesetting distinguishes between TST+ syntax and its semantics: compare respectively 
\begin{itemize*}
\item
$\tbot$ (syntax) and $\bot$ (a truth-value in $\{\bot,\top\}$); 
\item
$\timp$ and $\limp$; 
\item
$\tall$ and $\forall$; and 
\item
(deliberately making a prominent distinction, reflecting the design effort that goes into Definition~\ref{defn.model.tin}) $\tin$ and $\in$.
\end{itemize*}
It should therefore be clear when we are describing formal TST+ syntax, and when we are describing its meaning. 

We may use standard shorthands henceforth without comment, writing e.g. $\tneg\phi$ for $\phi\timp\tbot$; $\phi\tor\phi'$ for $(\tneg\phi)\timp\phi'$; $\phi\tand\phi'$ for $\tneg((\tneg\phi)\tor(\tneg\phi'))$; $\phi\tiff\phi'$ for $(\phi\timp\phi')\tand(\phi'\timp\phi)$; and $\texi a.\phi$ for $\tneg(\tall a.\tneg\phi)$.
\end{rmrk}

\begin{rmrk}
\begin{enumerate*}
\item
We will only be concerned with stratified predicates and terms henceforth, so all syntax is assumed stratified. 
We will not mention raw syntax again, except for Appendix~\ref{sect.aside} in which we show how non-normalising rewrites can appear if we relax the stratification condition.

Otherwise, if we write `predicate' then it is assumed to be stratified,
and if we write that $\phi$ is stratified then this is only for emphasis.
\item
Calling the sets of constants $\univi{i}$ is no coincidence: we have a particular value in mind for the $\univi{i}$ (Definition~\ref{defn.univi}).
However, for now $\univi{i}$ is just some arbitrary set of uninterpreted symbols for each $i\geq 0$.
\end{enumerate*}
\end{rmrk}

\begin{xmpl}[Examples of stratified and unstratified predicates]
Suppose $\level(b)=1$ and $\level(a)=2$ and $\level(c)=3$.
Then:
\begin{enumerate*}
\item
$\texi b.b\tin a$ is stratified. 
\item
$\tneg(a\tin a)$ and $\texi a.a\tin a$ and $c\tin a$ are unstratified.
\end{enumerate*}
More information can be found e.g. in~\cite{holmes:elestu,sep:connf}.
\end{xmpl}

\begin{nttn}[(Closed) predicates, free variables, and shift]
\label{nttn.plus}
\leavevmode
\begin{enumerate}
\item
Write $\tf{Pred}$ for the set of all predicates (the `stratified' here is understood now and henceforth).
\item\label{item.plus.closed}
Write $\f{fv}(\phi)$ for the free variables of $\phi$ as standard for first-order logic, and write $\consts(\phi)$ for the set of constants that appears in $\phi$.
Call a predicate $\phi\in\tf{Pred}$ \deffont{closed} when it has no free variables, and \deffont{open} when it does have free variables.
\item\label{item.upsome}
Given a closed predicate $\phi$ that mentions no constants, write $\phi\upsome{}$ (read as ``\deffont{shift $\phi$}'') for the predicate obtained by uniformly shifting all variable symbols up by one level.
So for instance:
$$
(\tall b.\texi a.b\tin a)\upsome{}=\tall b'.\texi a'.b'\tin a'
$$ 
where $\level(a')=\level(a)\plus 1$ and $\level(b')=\level(b')\plus 1$.
An inductive definition is not hard to write.
\end{enumerate}
\end{nttn}

\begin{defn}
\label{defn.nf.entailment}
\leavevmode
\begin{enumerate*}
\item
Let \deffont{TST entailment} be the usual derivation relation of first-order logic, along with the axiom \rulefont{Ext} 
in Figure~\ref{fig.TZT.axioms}.
\item
Let \deffont{TST+ entailment} be the derivation relation of TST, along with axioms
\begin{itemize*}
\item
\Cmp{\bar c, a,\phi} for every $a$ and $\phi$ with $\fv(\phi)\subseteq\{c_1,\dots,c_n,a\}$ ($\bar c$ abbreviates $c_1,\dots,c_n$), 
and 
\item
\rulefont{TA_\phi} in Figure~\ref{fig.TZT.axioms} for every closed $\phi$.
\end{itemize*}
\end{enumerate*}
\end{defn}

\begin{rmrk}
We add some comments on the axioms in Figure~\ref{fig.TZT.axioms} --- the footnotes annotate with technical notes which are not required for understanding on a first reading. 
\begin{enumerate}
\item
\rulefont{Ext} asserts \emph{extensionality}: that extensionally equal elements are equal.\footnote{Our version of TST+ does not include an explicit equality $\teq$.  As the reader may know, this loses no expressivity, because extensionality means we can treat $t\teq t'$ as a macro for $\tall a.a\tin t\tiff a\tin t'$.} 

Although we may call this an axiom, it is an axiom-scheme, with one (otherwise identical) axiom at each level.\footnote{In the presence of the typical ambiguity \rulefont{TA_\phi} axiom-scheme, it would suffice to include just one extensionality axiom, and let typical ambiguity do the rest.}
\item
\Cmp{\bar c, a,\phi} asserts \emph{comprehension}: for every $c_1$, \dots, $c_n$, the sets comprehension $\stxst{a}{\phi}$ \emph{``the set of $a$ such that $\phi$''} is indeed realised by some element.

Each $\phi$ gives rise to a distinct comprehension axiom, so this is an axiom-scheme in a somewhat stronger sense than \rulefont{Ext} is.\footnote{There are some things we could do to simplify this scheme, if we wish.  Notably: \emph{if} we know that our set of constants is equal to the underlying set of our intended denotation --- in this paper, that will be so and it is built in to Definition~\ref{defn.denotation} --- then we can save ourselves an ellipsis $\tall c_1\dots c_n$ and simplify to a scheme $\texi b.\tall a.a\tin b\tiff \phi$, where $\phi$ ranges over predicates (possibly mentioning constants) with $\fv(\phi)\subseteq\{a\}$.}
\item 
The $+$ in TST+ refers to its characteristic \deffont{typical ambiguity} axiom \rulefont{TA_\phi}, which asserts that validity (for closed predicates) is invariant under shifts of levels.

There is one \rulefont{TA_\phi} axiom for each choice of closed $\phi$ with no constants.\footnote{More technical details on how this works specifically in our proof, in Remark~\ref{rmrk.TA.consts}.}
\end{enumerate}
\end{rmrk}

\begin{figure}
$$
\begin{array}{l@{\qquad}l@{}l}
\rulefont{Ext}
&
\tall a,a'.(\tall b.b\tin a\tiff b\tin a') \timp \tall c.(a\tin c\tiff a'\tin c)
\\
\Cmp{\bar c, a,\phi}
&
\tall c_1\dots c_n.\texi b.\tall a.(a\tin b\tiff \phi) & \text{($\f{fv}(\phi)\subseteq\{c_1,\dots,c_n,a\}$)} 
\\
\rulefont{TA_\phi}
&
\phi\tiff \phi\upsome{} & (\phi\text{ closed, with no constants})
\end{array}
$$
\caption{TST+ axioms}
\label{fig.TZT.axioms}
\end{figure}

In summary: TST+ is a stratified first-order theory with extensionality, comprehension, and typical ambiguity. 
We now set about constructing a model for it.

\section{Hereditarily finitely supported sets}
\label{sect.hfs}

The point of this Section is to establish some background material on permutation sets and nominal sets.
The reader is welcome to read it now, or skip it at first and refer back as and when the material is used.
If there are two points we are aiming for in this Section, they are 
\begin{itemize*}
\item
to build the finitely supported powerset $\powersetfs(\ns X)$ (Definition~\ref{defn.pointwise.perm.action}(\ref{item.powersetfs})), and 
\item
to stack up a hierarchy of finitely supported powersets $\hfs{i}$ (Definition~\ref{defn.hfs}) and prove some basic nominal properties (Lemma~\ref{lemm.pi.pointwise}).
\end{itemize*}
$\hfs{}$ is essentially the same hierarchy as was used by Fraenkel and Mostowski to prove the independence of the Axiom of Choice from other other axioms of set theory~\cite{fraenkel:begdua,mostowski:unawo}.\footnote{Modulo some technical variations.  For instance, $\hfs{}$ in this paper is typed not cumulative, reflecting the stratification of TST (i.e. it is not the case that $\hfs{i}\subseteq\hfs{i\plus 1}$).} 
It is also the same sets universe on which nominal techniques are based~\cite{gabbay:newaas-jv,gabbay:fountl}.

\subsection{The basics: permutation sets and nominal sets}

\begin{defn}
\label{defn.set.of.atoms}
For the rest of this paper, fix a \deffont{set of atoms} $\atoms$ of size $\beth_{\omega}=\card\bigcup\powerset^i(\mathbb N)$, such that $\varnothing\not\in\atoms$ (this is a technical condition to avoid a pathology in Definition~\ref{defn.hfs}(\ref{item.hfs.0})).
\end{defn}

\begin{rmrk}
\label{rmrk.words.on.foundations}
Some words on how Definition~\ref{defn.set.of.atoms} relates to foundations:
\begin{enumerate}
\item
$\beth_\omega$ is the cardinality of the standard sets model of Higher-Order Arithmetic, obtained by starting from $\tf{Bool}$ and $\mathbb N$ and inductively forming all function-spaces.
$\beth_\omega$ is also the size of a model of Zermelo set theory.\footnote{This has the axioms of Infinity and Powerset, but unlike Zermelo\emph{-Fraenkel} set theory it lacks Replacement so that the function $i\in\mathbb N\longmapsto \powerset^i(\mathbb N)$ is a function-class but need not necessarily be a set.  Clearly, Zermelo set theory has much weaker consistency strength than Zermelo-Fraenkel set theory.}
\item
Atoms can be modelled using any set with cardinality $\beth_{\omega}$ such that $\varnothing\not\in\atoms$.
We do not require, and will not assume, any other structure to atoms (aside from being able to tell them apart).

If --- as argued for in~\cite{gabbay:equzfn} --- we work in ZFA (Zermelo-Fraenkel set theory with atoms), then elements in $\atoms$ could actually \emph{be} ZFA atoms.\footnote{This is nice because then atoms \emph{are atoms}, and from this follow by general principles some equivariance properties that in this paper we have to take a few minutes to by hand (such as Lemma~\ref{lemm.arrows.equivariant} or Lemma~\ref{lemm.supp.equivariant}).}
However, the maths that follows is self-contained and does not depend on a specific choice of $\atoms$: we just need 
the required cardinality.
\end{enumerate}
\end{rmrk}

\begin{defn}
\label{defn.perm}
\leavevmode
\begin{enumerate*}
\item
Write $\tf{Perm}$ for the set of permutations $\pi:\atoms\cong\atoms$ such that $\{\f{atom}\in\atoms \mid \pi(\f{atom})\neq\f{atom}\}$ is finite.
We call such a $\pi$ \deffont{finitely supported}.  

All permutations $\pi$ will be finitely supported henceforth. 
\item\label{item.perm.id}
Write $\id$ for the \deffont{identity permutation} such that $\id(\f{atom})=\f{atom}$ always.
\item\label{item.perm.circ}
If $\pi,\pi'\in\tf{Perm}$ then write $\pi\circ\pi'$ for their \deffont{composition}, defined by $(\pi\circ\pi')(\f{atom})=\pi(\pi'(\f{atom}))$.
\end{enumerate*}
\end{defn}

\begin{defn}
A \deffont{permutation set} $(\ns X,\act)$ consists of 
\begin{itemize*}
\item
a \deffont{carrier set} $\ns X$, along with 
\item
a \deffont{permutation action} $\act:\tf{Perm}\times\ns X\to \ns X$ that is a group action of $\tf{Perm}$ on $\ns X$ (so $\id\act x=x$ and $\pi\act(\pi'\act x)=(\pi\circ\pi')\act x$).
\end{itemize*}
\end{defn}

\begin{defn}
\label{defn.fix.S}
If $S\subseteq\atoms$ define $\fix(S)\subseteq\tf{Perm}$ by
$$
\fix(S)=\{\pi\in\tf{Perm} \mid \Forall{\f{atom}{\in}S} \mid \pi(\f{atom})=\f{atom}\} .
$$ 
\end{defn}

\begin{defn}
\label{defn.S.supports.x}
Suppose $(\ns X,\act)$ is a permutation set and $x\in\ns X$.
\begin{enumerate*}
\item\label{item.S.supports.x}
Say that $S\subseteq\atoms$ \deffont{supports} $x$, written $S\supports x$, when 
$$
\Forall{\pi{\in}\tf{Perm}}\pi\in\fix(S) \limp \pi\act x=x .
$$
In words: $S$ supports $x$ when if $\pi$ fixes $S$ pointwise --- meaning that $\pi\act\f{atom}=\f{atom}$ for every $\f{atom}\in S$ --- then $\pi$ fixes $x$.
\item\label{item.finitely.supported}
Call $x\in\ns X$ \deffont{finitely supported} when $S\supports x$ for some \emph{finite} $S\finsubseteq\atoms$.
\item
Call $(\ns X,\act)$ \deffont{finitely supported} when every $x\in\ns X$ is finitely supported.
\item
A \deffont{nominal set} is a finitely supported permutation set. 
\end{enumerate*}
\end{defn}

\begin{rmrk}
Note for experts: the terminology clash between the notion of finite support from Definition~\ref{defn.S.supports.x}(\ref{item.finitely.supported}) and the notion of a permutation $\pi$ being finitely supported used in Definition~\ref{defn.perm}, is only apparent:
$\tf{Perm}$ is naturally a nominal set under the \emph{conjugation action} $\pi'\act\pi=\pi'\circ\pi\circ(\pi')^\mone$ (which sends $\pi'(\f{atom})$ to $\pi\act\pi'(\f{atom})$) and with this action $\supp(\pi)=\{\f{atom}\in\atoms \mid \pi(\f{atom})\neq\f{atom}\}$.
\end{rmrk}

\begin{lemm}
Suppose $(\ns X,\act)$ is a permutation set.
Then:
\begin{enumerate*}
\item
If $x\in\ns X$ is finitely supported then $x$ has a unique least finite supporting set, which we write $\supp(x)$ and call the \deffont{support} of $x$.
\item
As a corollary, the permutation set $(\ns X,\act)$ is a nominal set precisely when $\supp$ is a total function from $\ns X$ to the finite powerset of $\atoms$.
\end{enumerate*}
\end{lemm}
\begin{proof}
The proof used here is modelled on \cite[Theorem~2.21]{gabbay:fountl}.
Suppose we have $A,B\finsubseteq\atoms$ such that $A\supports x$ and $B\supports x$; it will suffice to show that $(A\cap B)\supports x$; thus that if $\pi\in\fix(A\cap B)$ then $\pi\act x=x$.
Consider $\pi\in\fix(A\cap B)$ and let $\tau$ be self-inverse permutation that bijects $B\setminus A$ with some subset of $\atoms\setminus(A\cup B)$ and leave all other atoms fixed (so that $\tau=\tau^\mone$).
Write $\pi^\tau=\tau\circ\pi\circ\tau$ and note (using the group action) that $\pi\act x = \tau \act (\pi^\tau \act (\tau\act x))$.
Now $\tau\act x=x$ (because $\tau\in\fix(A)$), and $\pi^\tau\act x=x$ (because $\pi^\tau\in\fix(B)$).
Simplifying using these two equalities, it follows that 
$$
\pi\act x=\tau \act (\pi^\tau \act (\tau\act x))=x
$$ 
as required.

The corollary follows, since if an element has a unique least finite supporting set then it certainly has some finite supporting set, and if an element has some finite supporting set then by part~1 of this result it has a unique least such.
\end{proof}

\begin{defn}
\label{defn.pointwise.perm.action}
Suppose $(\ns X,\act)$ is a permutation set.
Then:
\begin{enumerate*}
\item
The powerset $\powerset(\ns X)$ is naturally a permutation set with the \deffont{pointwise} permutation action on $U\subseteq\ns X$ given by
$$
\pi\act U = \{\pi\act x \mid x\in U\} . 
$$ 
\item\label{item.powersetfs}
Write $\powersetfs(\ns X)\subseteq\powerset(\ns X)$ for the nominal set obtained by restricting to the finitely-supported subsets:
$$
\powersetfs(\ns X) = \{U\subseteq\ns X \mid \Exists{S\finsubseteq\atoms} S\supports U\} . 
$$ 
\end{enumerate*}
\end{defn}

\begin{lemm}
\label{lemm.supp.equivariant}
Suppose $(\ns X,\act)$ is a permutation set and $x\in\ns X$ and $\pi\in\tf{Perm}$.
Then:
\begin{enumerate*}
\item
If $S\finsubseteq\atoms$ then $S\supports x$ if and only if $\pi\act S\supports \pi\act x$.
\item 
$\supp(\pi\act x) = \pi\act\supp(x) = \{\pi(\f{atom}) \mid \f{atom}\in\supp(x)\}$. 

In words we say: $\supp$ is \emph{equivariant}, meaning that it commutes with the permutation action.
\item
If $x$ is finitely-supported then so is $\pi\act x$.
\item
As a corollary, $\powersetfs(\ns X)$ is a nominal set. 
\end{enumerate*}
\end{lemm}
\begin{proof}
By routine calculations from the group action.
\end{proof}

\subsection{A hierarchy of nominal sets}

We now know what a permutation set is, what a nominal set is, and how to build a finitely supported powerset.
We now iterate finitely supported powersets in a standard way, to build a hierarchy $\hfs{i}$ of hereditarily finitely supported sets.

\begin{defn}
\label{defn.hfs}
We define $\hfs{i}$ a \deffont{hereditarily finitely supported} hierarchy of nominal sets for $0\leq i\in\mathbb N$ as follows:
\begin{enumerate*}
\item\label{item.hfs.0}
$\hfs{0}=\atoms\cup\{\varnothing\}$ (we assumed $\varnothing\not\in\atoms$ in Definition~\ref{defn.set.of.atoms}) with the natural action given by
$$
\Forall{\f{atom}\in\atoms}\pi\act\f{atom}=\pi(\f{atom})
\quad\text{and}\quad
\pi\act\varnothing=\varnothing .
$$
\item
$\hfs{i\plus 1}=\powersetfs(\hfs{i})$, with 
its standard pointwise action $\pi\act X=\{\pi\act x\mid x\in X\}$ (Definition~\ref{defn.pointwise.perm.action}(1\&2)). 
\end{enumerate*} 
\end{defn}

\begin{rmrk}
Note for experts: We set $\hfs{0}$ to be $\atoms\cup\{\varnothing\}$ because including $\varnothing$ at level $0$ will be technically convenient later as a default value in clause~\ref{item.downarrow} of Definition~\ref{defn.uparrow.and.downarrow}.
However, this is also a principled choice: $\hfs{0}$ collects all the naturally extensionally empty elements of a universe-with-atoms.
\end{rmrk}

\begin{rmrk}
\label{rmrk.facts.about.hfs}
Unpacking definitions, if $\pi\in\tf{Perm}$ then each $\hfs{i}$ has a pointwise permutation action $\pi\act x$ such that: 
\begin{enumerate*}
\item
$\pi\act \f{atom}=\pi(\f{atom})$ for $\f{atom}\in\hfs{0}$, and $\pi\act\varnothing=\varnothing$.
\item
$\pi\act X=\{\pi\act x \mid x\in X\}$ for $X\in\hfs{i\plus 1}$.
\end{enumerate*}
We get finite support for every element at every stage by the construction, since:
\begin{itemize*}
\item
At level $0$ the reader can easily check that $\supp(\f{atom})=\{\f{atom}\}$ and $\supp(\varnothing)=\varnothing$.
\item
At level $i\plus 1$ we only admit $X\subseteq\hfs{i}$ into $\hfs{i\plus 1}$ when $X$ has finite support.
\end{itemize*}
\end{rmrk}

Lemma~\ref{lemm.pi.pointwise} gathers together some standard facts about the construction above:
\begin{lemm}
\label{lemm.pi.pointwise}
Suppose $i\geq 0$ and $x,x'\in\hfs{i}$ and $X\in\hfs{i\plus 1}$, and suppose $\pi,\pi'\in\tf{Perm}$ are permutations.
Then: 
\begin{enumerate*}
\item\label{item.pi.pointwise.X}
$\pi\act X=\{\pi\act x\mid x\in X\}$. 
\item\label{item.pi.pointwise}
$(\pi\act x)\in (\pi\act X)$ if and only if $x\in X$.
\item\label{item.pi.group}
The action is indeed a group action: $\id\act x=x$ and $\pi\act (\pi'\act x)=(\pi\circ\pi')\act x$ (Definition~\ref{defn.perm}(\ref{item.perm.id}\&\ref{item.perm.circ})). 
\item\label{item.pi.bij}
As a corollary, the action $x\mapsto \pi\act x$ is bijective on $\hfs{i}$.
\item\label{item.supp.supports}
If $\pi\in\fix(\supp(x))$ then $\pi\act x=x$ ($\fix$ is defined in Definition~\ref{defn.fix.S}).
\end{enumerate*}
\end{lemm}
\begin{proof}
We consider each part in turn:
\begin{enumerate}
\item
This is just Definition~\ref{defn.pointwise.perm.action} and Remark~\ref{rmrk.facts.about.hfs}(2).
\item
From part~\ref{item.pi.pointwise.X} of this result.
\item
By a routine induction on $\hfs{}$.
The base case is from Remark~\ref{rmrk.facts.about.hfs}(1) and for the induction step we can use Remark~\ref{rmrk.facts.about.hfs}(2).
\item
From part~\ref{item.pi.group} of this result, since $x\mapsto\pi^\mone\act x$ is inverse to $x\mapsto \pi\act x$. 
\item
From Definition~\ref{defn.S.supports.x}(\ref{item.S.supports.x}), since by construction $\supp(x)\supports x$.
\qedhere\end{enumerate}
\end{proof}

\begin{xmpl}\leavevmode
\label{xmpl.hfs.examples}
We compute some values for $\hfs{}$ from Definition~\ref{defn.hfs}:
\begin{enumerate*}
\item
$\hfs{0}$ is just $\atoms\cup\{\varnothing\}$.
\item
$\hfs{1}$ is the set of finite and cofinite subsets of $\hfs{0}$.
\item
$\hfs{2}$ is \emph{not} just the set of finite and cofinite subsets of $\hfs{1}$.
For example, $\hfs{2}$ contains $\mathbf i_2=\{x\subseteq\hfs{1} \mid \card x=i\}$ for every $i\in\mathbb N$.
Note that $\supp(\mathbf i_2)=\varnothing$ for every $i$ (so it has empty support, even though its elements do not). 
\item
$\mathbb N_2=\{\mathbf i_2 \mid i\in\mathbb N\}\in\hfs{3}$.
Indeed, every subset $X\subseteq \mathbb N_2$ is also in $\hfs{3}$ (the reader can check that every such $X$ has empty support).
\item
$\hfs{3}$ also contains $\mathbf i_3=\{x\subseteq\hfs{2} \mid \card x=i\}$ for every $i\in\mathbb N$, and plenty of other structure besides, such as (for example) $\{\{\{a\}\}\}$ for every $a\in\atoms$.
So this is already quite a complicated set.
\item
\dots and so on.
\end{enumerate*}
\end{xmpl}

\section{Two ways to map between $\hfs{\lowercasei}$ and $\hfs{\lowercasei\plus 1}$}
\label{sect.map.between}

In order to model the \emph{typical ambiguity} property of TST+ (Definition~\ref{defn.nf.entailment}; Proposition~\ref{prop.TA}), it will be useful to be able to map between $\hfs{i}$ and $\hfs{i\plus 1}$.
In fact, we will use three such maps:
\begin{enumerate*}
\item
$\uparrow$ and $\downarrow$, defined in Definition~\ref{defn.uparrow.and.downarrow}, are raising and lowering functions based on general facts about sets structure.
They are fairly straightforward: $\uparrow$ will be injective and $\downarrow$ will be surjective.
\item
$f$ is a bijection between $\hfs{i}$ and $\hfs{i\plus 1}$ obtained by a bespoke nominal cardinality argument.
The reader might assume that this could not exist, since $\hfs{i\plus 1}$ is a powerset of $\hfs{i}$, but it is a \emph{finitely supported} powerset, which makes all the difference (see Remark~\ref{rmrk.external.f}).
\end{enumerate*}

\subsection{Raising $\uparrow:\hfs{i}\to\hfs{i\plus 1}$ and lowering $\downarrow:\hfs{i\plus 1}\to\hfs{i}$} 

\begin{defn}
\label{defn.uparrow.and.downarrow}
For each $i\geq 0$ we define 
\begin{itemize*}
\item
a \deffont{raising} function $x\in\hfs{i}\longmapsto x{\uparrow}\in\hfs{i\plus 1}$, and
\item 
a \deffont{lowering} function $x\in\hfs{i\plus 1}\longmapsto x{\downarrow}\in\hfs{i}$
\end{itemize*}
by mutual induction as follows: 
\begin{enumerate}
\item\label{item.uparrow.0}
If $x\in\hfs{0}$ then
$$
x{\uparrow}=\{x\} 
$$
\item\label{item.downarrow}
If $X\in\hfs{1}$ then
$$
X{\downarrow}=
\begin{cases}
x & \text{if $X=\{x\}$ for some $x\in\hfs{0}$}
\\
\varnothing & \text{otherwise}
\end{cases}
$$
\item\label{item.uparrow}
If $x\in\hfs{i\plus 1}$ for $i\geq 0$ then 
$$
x{\uparrow}=\{y\in\hfs{i\plus 1} \mid y{\downarrow}\in x\} 
$$
\item
If $X\in\hfs{i\plus 1}$ for $i\geq 0$ then 
$$
X{\downarrow}=\{y\in\hfs{i} \mid y{\uparrow}\in X\}  
$$
\end{enumerate}
We may indicate multiple applications of $\uparrow$ and $\downarrow$ with subscripts, e.g. $x{\uparrow}^2= (x{\uparrow}){\uparrow}$.
\end{defn}

\begin{rmrk}
We give examples of $\uparrow$ and $\downarrow$ in action in Example~\ref{xmpl.updown.examples}.

Lemmas~\ref{lemm.arrows.equivariant}, \ref{lemm.arrows.adjoint}, and~\ref{lemm.x.up} develop properties of $\uparrow$ and $\downarrow$ by standard sets manipulations.

It certainly looks like ${\uparrow}:\hfs{i}\to\hfs{i\plus 1}$ and ${\downarrow}:\hfs{i\plus 1}\to\hfs{i}$ in Definition~\ref{defn.uparrow.and.downarrow} but this does need to be checked: this happens in Lemma~\ref{lemm.arrows.equivariant}(2).
\end{rmrk}

\begin{lemm}
\label{lemm.arrows.equivariant}
\leavevmode
\begin{enumerate*}
\item
Raising $\uparrow$ and lowering $\downarrow$ are \emph{equivariant}, by which we mean that if $\pi\in\tf{Perm}$ then 
$$
\pi\act (x{\uparrow})=(\pi\act x){\uparrow}
\quad\text{and}\quad
\pi\act (x{\downarrow})=(\pi\act x){\downarrow}.
$$
\item
As a corollary, raising $\uparrow:\hfs{i}\to\hfs{i\plus 1}$ and lowering $\downarrow:\hfs{i\plus 1}\to\hfs{i}$ do indeed map between hereditarily finitely supported powersets. 
\end{enumerate*}
\end{lemm}
\begin{proof}
By routine inductions.
We sketch the argument:
\begin{enumerate}
\item
We work by induction on $i$:
\begin{itemize*}
\item
For the case of $x{\uparrow}$ for $x\in\hfs{0}$ we unfold definitions and use the pointwise action: 
$$
\pi\act(x{\uparrow})=\pi\act\{x\}=\{\pi\act x\}=(\pi\act x){\uparrow}.
$$
\item
For the case of $X{\uparrow}$ for $X\in\hfs{i\plus 1}$ where $i\geq 0$, we just chase definitions and use the pointwise action and the inductive hypothesis.
We sketch the reasoning: 
\begin{multline*}
\pi\act(X{\uparrow})\stackrel{\text{def}}=\pi\act\{x\in\hfs{i}\mid x{\downarrow}\in X\}\stackrel{\text{pointwise act.}}=\{\pi\act x\mid x{\downarrow}\in X\}
\stackrel{\text{perm. act. bijective}}= 
\\
\{x\mid (\pi^\mone\act x){\downarrow}\in X\}\stackrel{\text{ind. hyp. ${\downarrow}$}}=\{x\mid \pi^\mone\act(x{\downarrow})\in X\}\stackrel{\text{L.\ref{lemm.pi.pointwise}(\ref{item.pi.pointwise})}}=\{x\mid x{\downarrow}\in\pi\act X\} \stackrel{\text{def}}= (\pi\act X){\downarrow}
\end{multline*}
\end{itemize*}
The argument for $\downarrow$ is similar and no harder.
\item
It suffices to show that $\supp(x{\uparrow})\subseteq\supp(x)$ and $\supp(x{\downarrow})\subseteq\supp(x)$.
This is easy from part~1 of this result, since if $\pi\in\fix(\supp(x))$ then $\pi\act (x{\uparrow})=(\pi\act x){\uparrow}=x{\uparrow}$, and similarly for $x{\downarrow}$.
\qedhere\end{enumerate}
\end{proof}

\begin{lemm}
\label{lemm.arrows.adjoint}
$\uparrow$ and $\downarrow$ are \deffont{adjoint}, by which we mean that for $i\geq 0$:
\begin{enumerate*}
\item
If $y,x\in\hfs{i\plus 1}$ then
$y{\downarrow}\in x$ if and only if $y\in x{\uparrow}$.
\item
If $y\in\hfs{i}$ and $X\in\hfs{i\plus 2}$ then
$y{\uparrow}\in X$ if and only if $y\in X{\downarrow}$.
\end{enumerate*}
\end{lemm}
\begin{proof}
Directly from the construction in Definition~\ref{defn.uparrow.and.downarrow}.
\end{proof}

\begin{lemm}
\label{lemm.x.up}
Suppose $i\geq 0$.
Then:
\begin{enumerate*}
\item\label{x.up.down.x}
If $x\in\hfs{i}$ then $x{\uparrow}{\downarrow}=x$.
\item\label{item.uparrow.injective}
$\uparrow:\hfs{i}\to\hfs{i\plus 1}$ is injective: 
if $x,x'\in\hfs{i}$ and $x{\uparrow}=x'{\uparrow}$ then $x=x'$.
\item\label{item.downarrow.surjective}
$\downarrow:\hfs{i\plus 1}\to\hfs{i}$ is surjective: 
every $x\in\hfs{i}$ is equal to $X{\downarrow}$ for some $X\in\hfs{i\plus 1}$. 
\end{enumerate*}
\end{lemm}
\begin{proof}
We prove part~\ref{x.up.down.x} by an easy induction on $i$:
\begin{itemize*}
\item
The case that $i=0$ is just from clauses~1 and~2 of Definition~\ref{defn.uparrow.and.downarrow} since $\{x\}{\downarrow}=x$.
\item
For the case of $i\geq 1$ we note using Lemma~\ref{lemm.arrows.adjoint} and the inductive hypothesis that 
$$
y\in x{\uparrow}{\downarrow} \liff y{\uparrow}\in x{\uparrow} \liff y{\uparrow}{\downarrow}\in x \stackrel{ind. hyp.}\liff y\in x.
$$
It follows that $x{\uparrow}{\downarrow}=x$.
\end{itemize*}
Injectivity of ${\uparrow}$ follows, since if $x{\uparrow}=x'{\uparrow}$ then also $x=x{\uparrow}{\downarrow}=x'{\uparrow}{\downarrow}=x'$ (i.e. ${\uparrow}$ has a right inverse).
Surjectivity of ${\downarrow}$ follows, since $x=(x{\uparrow}){\downarrow}$ (i.e. ${\downarrow}$ has a left inverse). 
\end{proof}

\begin{xmpl}
\label{xmpl.updown.examples}
All we need about $\uparrow$ and $\downarrow$ from Definition~\ref{defn.uparrow.and.downarrow} is what is stated in the lemmas above: this is a sequence of pairs of adjoint maps that are equivariant and injective going upwards and surjective going downwards.
Beyond possessing these properties, we will not need to care what values $\uparrow$ and $\downarrow$ actually \emph{produce}. 
That said, there is no harm in trying out a few values, and the adjoint property from Lemma~\ref{lemm.arrows.adjoint} makes these relatively straightforward to compute.

So suppose $a\in\mathbb A\subseteq\hfs{0}$ and $x\in\hfs{1}$. 
Then:
\begin{enumerate}
\item
$x\in a{\uparrow}{\uparrow}$ if and only if $x{\downarrow}\in\{a\}$ if and only if $x{\downarrow}=a$ if and only if $x=\{a\}$.
Thus, $a{\uparrow}{\uparrow}=\{\{a\}\}$.
\item
$x\in \varnothing{\uparrow}{\uparrow}$ if and only if $x{\downarrow}\in\{\varnothing\}$ if and only if $x{\downarrow}=\varnothing$ if and only if $x$ is not a singleton set of an element in $\hfs{0}$.
Thus, $\varnothing{\uparrow}{\uparrow}=\{x\in\hfs{1} \mid \card x \neq 1\}$.
\item
Suppose $X\in\hfs{3}$.
Then $a\in X{\downarrow}{\downarrow}$ if and only if $\{a\}\in X{\downarrow}$ if and only if $\{\{a\}\}\in X$, and $\varnothing\in X{\downarrow}{\downarrow}$ if and only if $\{x\in\hfs{1}\mid \card x\neq 1\}\in X$.
Thus, $X{\downarrow}{\downarrow}=\{a\in\atoms \mid \{\{a\}\}\in X\}$ or $X{\downarrow}{\downarrow}=\{a\in\atoms\mid\{\{a\}\}\in X\}\cup\{\varnothing\}$, depending on whether $\{x\in\hfs{1}\mid \card x\neq 1\}\in X$.
\end{enumerate}
\end{xmpl}

\subsection{The bijection $f:\hfs{i}\cong\hfs{i\plus 1}$}

\begin{rmrk}
\label{rmrk.sketch.proof}
The construction of $f:\hfs{i}\cong\hfs{i\plus 1}$ is a bit fiddly, in the way of many counting arguments, so we take a moment to suggest at a high level how the bijection works. 

The finite support condition on elements $X\in\hfs{i\plus 1}$ in Definition~\ref{defn.hfs}(2) does restrict cardinalities: e.g. as we noted in Example~\ref{xmpl.hfs.examples}, $\hfs{1}$ is just the finite and cofinite subsets of $\hfs{0}$, and this has the same cardinality as $\hfs{0}$, namely $\beth_\omega$ --- whereas the full powerset $\powerset(\hfs{0})$ has cardinality $\beth_{\omega\plus 1}$.

Looking more closely, we see that an element in $\hfs{i}$ can be represented as what we might call a \emph{colourless sketch} along with information about how to \emph{colour in} that sketch with particular atoms.

As we take finitely supported powersets the cardinality of the set of possible sketches increases --- but if the set of atoms (`colours') is large enough, there will still be more colours than sketches and the size of the set of colours will dominate.

Or to put it another way: we made the set of atoms so large in Definition~\ref{defn.set.of.atoms} that all of the detectable cardinality of $\hfs{i}$ comes from choices of atoms, and none of it comes from actual sets structure.

This needs to be made formal and checked of course --- but the basic ideas are there, and indeed we can see these ideas quite faithfully represented in the maths: see Definitions~\ref{defn.X|S} (for the `sketch' structure) and~\ref{defn.colouring.in} (for the `colouring in').

The culminating result of this Section is Corollary~\ref{corr.exists.f}, where we observe that a (non-equivariant) bijection $f:\hfs{i}\to\hfs{i\plus 1}$ exists.

Note that the \emph{only} fact we need about $f$ for the later proofs is that it exists --- we do not need it to be equivariant (it is not), and we do not need it to do anything more than just biject sets.
Thus, the reader who is sufficiently convinced by the argument above can safely skip to Corollary~\ref{corr.exists.f} and proceed to the next Section.
\end{rmrk}

\begin{defn}
\label{defn.B.choice}
For this Subsection, fix:
\begin{enumerate*}
\item
A countable subset $\btoms\subseteq\atoms$ (so $\card\btoms=\beth_0=\card\mathbb N$).
\item
A single atom $\f{notinB}\in\atoms\setminus\btoms$.
\end{enumerate*}
\end{defn}

Definition~\ref{defn.X|S} is parameterised over our choices of $\mathbb B$ and $\f{notinB}\in\atoms\setminus\mathbb B$ in Definition~\ref{defn.B.choice}.
We never need to vary this choice, so $\mathbb B$ and $\f{notinB}$ are constant henceforth:
\begin{defn}
\label{defn.X|S}
Suppose $X\in\hfs{i}$ for $i\geq 0$.
Define $X|_\btoms$ --- call this the \deffont{sketch} of $X$, for reasons discussed in Remark~\ref{rmrk.sketch.proof}
--- by induction on $i$ as follows:\footnote{$X|_\btoms$ is an instance of a general family of constructions.  For example, compare the $X|_\btoms$ defined here with the related $X|_{\#a}$ from \cite[Definition~8.8]{gabbay:fountl}.  In a nominal context, these constructions can be quite versatile.}
\begin{enumerate*}
\item
If $X\in\hfs{0}$ then 
$$
X|_\btoms = 
\begin{cases}
X & X\in\btoms 
\\
\f{notinB} & X\not\in\btoms 
\end{cases}
$$
\item
If $X\in\hfs{i}$ for $i\geq 1$ then 
$$
X|_\btoms=\{x|_\btoms \mid x\in X,\ \supp(x)\subseteq \btoms\}.
$$
\qedhere\end{enumerate*}
\end{defn}

\begin{xmpl}
\label{xmpl.colour}
Suppose $a,b\in\atoms$ and $a\not\in\mathbb B$ and $b\in\mathbb B$.
Then the reader can check of Definition~\ref{defn.X|S} that:
\begin{enumerate*}
\item
$a|_{\mathbb B}=\f{notinB}=\varnothing|_{\mathbb B}$ and $b|_{\mathbb B}=b$.
\item
$\atoms\in\hfs{1}$ and $\atoms|_{\mathbb B} = \mathbb B$.
\item
$\{\atoms\}\in\hfs{2}$ and $\{\atoms\}|_{\mathbb B} = \{\mathbb B\}$.
\item
$\{a,b\}\in\hfs{1}$ and $\{a,b\}|_{\mathbb B}=\{b\}$.
\item
$\atoms\setminus\{a,b\}\in\hfs{1}$ and $(\atoms\setminus\{a,b\})|_{\mathbb B}=\mathbb B\setminus \{b\}$. 
\item
$\{\varnothing,a,b\}\in\hfs{1}$ and $\{\varnothing,a,b\}|_{\mathbb B}=\{\f{notinB},b\}$.
\item
$\{\{\varnothing\},\{a\},\{b\},\atoms{\setminus}\{a\},\atoms{\setminus}\{b\}\}\in\hfs{2}$ and 
$$
\{\{\varnothing\},\{a\},\{b\},\atoms{\setminus}\{a\},\atoms{\setminus}\{b\}\}|_{\mathbb B}=\{\{\f{notinB}\},\{b\},\mathbb B{\setminus}\{b\}\}.
$$
\end{enumerate*}
\end{xmpl}

\begin{lemm}
\label{lemm.X.X'}
Suppose $x,x'\in\hfs{0}$ and $\supp(x),\supp(x')\subseteq\btoms$.
Then
$$
x=x'
\quad\text{if and only if}\quad x|_\btoms=x'|_\btoms.
$$
\end{lemm}
\begin{proof}
It is immediate that $x=x'$ implies $x|_\btoms=x'|_\btoms$. 

Suppose $x\neq x'$.
We consider cases:
\begin{itemize}
\item
$x=\varnothing$ and $x'=\varnothing$ is impossible because we assumed $x\neq x'$.
\item
$x\in\atoms\setminus\btoms$ or $x'\in\atoms\setminus\btoms$ is impossible because we assumed $\supp(x),\supp(x')\subseteq\btoms$ and in this case (as noted in Remark~\ref{rmrk.facts.about.hfs}) $\supp(x)=\{x\}$ and $\supp(x')=\{x'\}$.
\item
Suppose $x=\varnothing$ and $x'\in\btoms$.
Then $x|_\btoms=\f{notinB}$ and $x'|_\btoms=x'$.
Thus $x'\neq\f{notinB}$ just because $\f{notinB}$ is not in $\mathbb B$. 
\item
The case that $x'=\varnothing$ and $x\in\btoms$ is symmetric.
\item
Suppose $x,x'\in\btoms$.
Then $x|_\btoms=x$ and $x'|_\btoms=x'$.
Thus $x|_\btoms\neq x'|_\btoms$.
\qedhere\end{itemize} 
\end{proof}

\begin{lemm}
\label{lemm.X|S.facts}
Suppose $i\geq 0$.
Then:
\begin{enumerate*}
\item
Suppose $X,X'\in\hfs{i}$ and $\supp(X),\supp(X')\subseteq\btoms$.
Then 
$$
X=X'
\quad\text{if and only if}\quad X|_\btoms=X'|_\btoms.
$$
\item
Suppose $x\in\hfs{i}$ and $X\in\hfs{i\plus 1}$ and $\supp(x),\supp(X)\subseteq\btoms$.
Then 
$$
x\in X
\quad\text{if and only if}\quad
x|_\btoms\in X|_\btoms.
$$
\end{enumerate*}
\end{lemm}
\begin{proof}
We prove both parts by a simultaneous induction on $i$.
\begin{enumerate}
\item
It is immediate that $X=X'$ implies $X|_\btoms=X'|_\btoms$. 

Now suppose $X\neq X'\in\hfs{i}$.
If $i=0$ then we use Lemma~\ref{lemm.X.X'} above.

So suppose $i\geq 1$.
By sets extensionality there exists $x\in\hfs{i\minus 1}$ such that $x\in X\land x\not\in X'$ or $x\not\in X\land x\in X'$; we consider the former case (the latter case is symmetric).

Recall that we assumed that $\supp(X),\supp(X')\subseteq\btoms$.
Let $\tau\in\fix(\supp(X)\cup\supp(X'))$ be a self-inverse permutation that bijects $\supp(x)\setminus(\supp(X)\cup\supp(X'))$ with some finite subset of $\btoms\setminus(\supp(X)\cup\supp(X'))$ (so $\tau=\tau^\mone$).
If $\supp(x)\subseteq\btoms$ already, then we can just set $\tau=\id$.
This $\tau$ exists, because $\btoms$ is countable and $\supp(x)$, $\supp(X)$, and $\supp(X')$ are all finite.

By Lemma~\ref{lemm.supp.equivariant} $\supp(\tau\act x)=\tau\act\supp(x)\subseteq\btoms$. 
By Lemma~\ref{lemm.pi.pointwise}(\ref{item.pi.pointwise}) $x\in X$ if and only if $\tau\act X\in\tau\act X$ and similarly for $X'$, and by Lemma~\ref{lemm.pi.pointwise}(\ref{item.supp.supports}) (since $\tau\in\fix(\supp(X)\cup\supp(X'))$) $\tau\act X=X$ and $\tau\act X'=X'$.

By part~2 of the inductive hypothesis at level $i\minus 1$ we have that $(\tau\act x)|_\btoms\in X|_\btoms$ and $(\tau\act x)|_\btoms\not\in X'|_\btoms$. 
Thus, $X|_\btoms\neq X'|_\btoms$. 
\item
$x\in X$ implies $x|_\btoms\in X|_\btoms$ follows from Definition~\ref{defn.X|S}.

Now suppose $x\not\in X\in\hfs{i\plus 1}$ and $x|_\btoms\in X|_\btoms$.
This means that there exists $x'\in X$ such that $\supp(x')\subseteq\btoms$ and $x'|_\btoms = x|_\btoms$.
We assumed $\supp(x)\subseteq\btoms$ so by part~1 of the inductive hypothesis at level $i$, it follows that $x=x'$, which is a contradiction since we assumed $x\not\in X$.
\qedhere\end{enumerate}
\end{proof}

\begin{defn}
For each $i\geq 1$ we define
$$
\hfs{i}|_\btoms = \{x|_\btoms \mid x\in\hfs{i},\ \supp(x)\subseteq\btoms\} .
$$
This is indeed equal to $\hfs{i}|_\btoms$ in the sense of Definition~\ref{defn.X|S}(2), for $\hfs{i}$ considered as an element of $\hfs{i\plus 1}$.
\end{defn}

\begin{defn}
\label{defn.colouring.in}
Suppose $i\geq 0$ and $(\pi,x')\in\tf{Perm}\times\hfs{i}|_\btoms$.
Then define $\int_\pi x'$ --- we might call this the \deffont{colouring in} of $x'$ by $\pi$ --- as follows:
\begin{itemize*}
\item
$\int_\pi x'$ is the (by Lemma~\ref{lemm.X|S.facts}(2) unique, if it exists) $x\in\hfs{i}$ such that 
$$
\supp(\pi^\mone\act x)\subseteq\btoms
\quad\text{and}\quad 
(\pi^\mone\act x)|_\btoms = x'.
$$
\item
$\int_\pi x'=\varnothing\in\hfs{i}$ if no such $x$ exists.
\end{itemize*}
\end{defn}

\begin{prop}
\label{prop.pi.act.surjection}
Suppose $i\geq 1$. 
Then 
$$
(\pi,x')\in\tf{Perm}\times\hfs{i}|_\btoms\longmapsto\textstyle\int_\pi x'\in\hfs{i}
$$ 
is a surjection from $\tf{Perm}\times\hfs{i}|_\btoms$ to $\hfs{i}$. 
\end{prop}
\begin{proof}
For $x\in\hfs{i}$, let $\pi\in\tf{Perm}$ be any permutation that maps $\supp(x)$ into $\btoms$; this exists because $\supp(x)$ is finite and $\btoms$ is countably infinite and by Lemma~\ref{lemm.supp.equivariant} $\supp(\pi\act x)=\pi\act\supp(x)=\{\pi(\f{atom})\mid \f{atom}\in\supp(x)\}$.

Then $(\pi^\mone,(\pi\act x)|_\btoms)\in\tf{Perm}\times\hfs{i}|_\btoms$ and this maps to $x$. 
\end{proof}

\begin{lemm}
\label{lemm.card.limit}
Suppose $i\geq 1$.
Then $\card(\tf{Perm}\times\hfs{i}|_\btoms)=\beth_\omega=\card\atoms$.
\end{lemm}
\begin{proof}
$\tf{Perm}$ from Definition~\ref{defn.perm}(1) is a set of \emph{finitely supported} permutations of $\atoms$, and it is a fact that it follows that $\card\tf{Perm}=\card\atoms$.
$\hfs{i}|_\btoms$ is a subset of a finitely iterated powerset of a \emph{countable} set (namely $\btoms\cup\{\f{notinB}\}$), therefore it has cardinality $\card\powerset^i(\mathbb N)$ and this is strictly less than $\card\atoms=\beth_\omega=\card\bigcup_{i\in\mathbb N}\powerset^i(\mathbb N)$. 
\end{proof}

\begin{corr}
\label{corr.exists.f}
\leavevmode
\begin{enumerate*}
\item
$\hfs{i}$ bijects with $\atoms$, for every $i\geq 0$.
\item
As a corollary, for every $i\geq 0$ there exists a bijection 
$$
f_i:\hfs{i}\to\hfs{i\plus 1}.
$$
\end{enumerate*}
We will usually drop the subscript and write $f_i$ just as $f$, since the level will always be unimportant or understood.\footnote{Note that there is no requirement that $f$ be equivariant or that it preserve or reflect $\in$-structure, and in general $f$ does neither: it need not be the case that $f(\pi\act x)=\pi\act f(x)$, nor that $y\in x$ if and only if $f(y)\in f(x)$.}
\end{corr}
\begin{proof}
It would suffice to show that $\card\hfs{i}=\card\mathbb A$.
If $i=0$ then $\hfs{i}=\atoms\cup\{\varnothing\}$ and the result follows.
So suppose $i\geq 1$.
We prove two inequalities:
\begin{itemize}
\item
$\card\hfs{i}\geq\card\atoms$ because there is an obvious injection from $\atoms$ into $\hfs{i}$ given by taking iterated singletons; e.g. $\f{atom}\in\atoms$ maps to $\{\{\f{atom}\}\}\in\hfs{2}$.
\item
Conversely, by Proposition~\ref{prop.pi.act.surjection} $\tf{Perm}\times\hfs{i}|_\btoms$ surjects onto $\hfs{i}$, and by Lemma~\ref{lemm.card.limit} $\card(\tf{Perm}\times\hfs{i}|_\btoms)=\card\atoms$.
Thus also $\card\hfs{i}\leq\card\atoms$, and we are done.
\qedhere\end{itemize}
\end{proof}

\begin{rmrk}
\label{rmrk.external.f}
How can $f$ exist, given that it is bijecting a set with its (finitely supported) powerset?

We avoid a G\"odel diagonalisation argument because $f$ itself is not finitely supported --- which, for functions, means that there exists no finite $S\finsubseteq\atoms$ such that if $\pi\in\fix(S)$ (Definition~\ref{defn.fix.S}) then $\pi\act(f(x))=f(\pi\act x)$ for every $x$.

The point is that the bijection is external to $\hfs{}$: it cannot be represented as a graph within the $\hfs{}$ hierarchy.\footnote{Model theorists will be very well familiar with the distinction between classes and sets, and a function that exists outside a model and one that can be represented inside it.  This is precisely what is happening here.}
\end{rmrk}

\section{The model $\univi{}$ and denotations $\model{\phi}$ in that model}
\label{sect.the.model}

\subsection{Some intuitions}
\label{subsect.some.intuitions}

We now have all the parts we need to construct a universe $\univi{}=(\univi{i} \mid i\in\mathbb N)$ in which to interpret predicates of TST+.
Before we put these parts together, we give some high-level intuition of the construction.

$\univi{}$ is parameterised over a fixed but arbitrary number $n\in\mathbb N$ (Definition~\ref{defn.n}).
According to this parameter, $\univi{}$ (Definition~\ref{defn.univi}) splits into two halves 
$$
(\univi{i}\mid 0\leq i\leq n)
\quad\text{and}\quad
(\univi{n\plus i} \mid 0\leq i) ,
$$
as shown in Figure~\ref{fig.simple.diagram}. 
These halves overlap at $\univi{n}$; most of the action in the model happens there. 
\begin{itemize*}
\item
$\univi{i}$ is defined to be $\hfs{i}=\powersetfs^i(\atoms\cup\{\varnothing\})$ for $0\leq i\leq n$. 

Thus, the universes $\univi{i}$ at levels $i\leq n$ are a hierarchy 
of increasingly complex hereditarily finite sets (examples in Example~\ref{xmpl.hfs.examples}).

By Corollary~\ref{corr.exists.f} $\hfs{i}\cong\hfs{i\plus 1}$ for $0\leq i\leq n\minus 1$ using the bijection $f$, so these sets all have cardinality $\beth_\omega$.
However, their $\in$-structure increases in complexity with each level, up to $n$ ($f$ does not respect support or $\in$-structure; it is just a raw bijection between underlying sets).
\item
Each $\univi{n\plus i}$ for $i\geq 1$ is by construction a subset of $\univi{n\plus i\minus 1}$,
but bijections are given by $\uparrow$ and $\downarrow$ from Definition~\ref{defn.uparrow.and.downarrow}, for which it follows from adjointness (Lemma~\ref{lemm.arrows.equivariant}) and equivariance (Lemma~\ref{lemm.arrows.adjoint}) that they exhibit $\univi{n\plus i}$ as identical to $\univi{n}$ in all but name, for every $0\leq i$. 
\end{itemize*}

\begin{figure}
$$
\hspace{-1em}
\begin{array}{r@{}lcr@{}lcr@{}lcr@{}lcr@{}lc}
\univi{0}{=}&\atoms{\cup}\{\varnothing\}
&\stackrel{f}{\cong}&
\univi{1}{=}&\powersetfs(\univi{0})
&\stackrel{f}{\cong}& 
\univi{2}{=}&\powersetfs^2(\univi{0})
&\stackrel{f}{\cong}& 
\dots&
&\stackrel{f}{\cong}& 
\univi{n}{=}&\powersetfs^n(\univi{0})
\\[2ex]
\univi{n}{=}&\hfs{n}
&\stackrel{\uparrow}{\cong}&
\univi{n\plus 1}{=}&\{x{\uparrow}\mid x{\in}\univi{n}\}
&\stackrel{\uparrow}{\cong}& 
\univi{n\plus 2}{=}&\{x{\uparrow}^2\mid x{\in}\univi{n}\}
&\stackrel{\uparrow}{\cong}&
\dots&
\end{array}
$$
\caption{A simple map of $\univi{}$}
\label{fig.simple.diagram}
\end{figure}

As we work through the proofs, we will see that most of the really interesting structure of our model is determined at $\univi{n}$.
This includes the interpretation of $y\tin x$ in Definition~\ref{defn.model.tin}, and the `adjusted permutation action' and `adjusted support' from Definitions~\ref{defn.permutation.action} and~\ref{defn.modified.supp}.

Though the definitions may seem a little complex at first, all they do is shift an element $x\in\univi{i}$ to present it at $\univi{n}$ (up using $f$ if $i<n$, and down using $\downarrow$ if $i>n$) and then do whatever relevant operation needs done, at level $n$.

We still need the universes $\univi{i}$ for $i\neq n$, in order to build up to $\univi{n}$ and in order to interpret the stratified language of TST+ --- but in the background everything gets shifted to $\univi{n}$, processed there, and (if required) shifted back. 

This suffices to interpret TST+ derivations 
provided they are less complex than $n$, where complexity is made formal by a measure that we call \emph{extent}; see Notation~\ref{nttn.extent}.
For now, we just need to know that the extent of a derivation can be determined from its syntax and every derivation has finite extent.
Because $n$ is fixed but arbitrarily large, it follows by compactness of first-order logic that TST+ is consistent.\footnote{In Subsection~\ref{subsect.elaborating.the.model} we indicate how a model --- actually two quite large general classes of models --- capable of interpreting all of TST+ can easily be created from this using standard model-theoretic constructions (such as ultraproducts).}

\subsection{We construct our universe $\univi{}$}

\begin{defn}
\label{defn.n}
Make a fixed but arbitrary choice of number $n\geq 1$.
This will be a parameter in the model we build in Definition~\ref{defn.univi}, and will be fixed throughout the maths that follows.
The proofs use $n$ substantively just once, in Proposition~\ref{prop.comprehension}.
\end{defn}

\begin{defn}
\label{defn.univi}
We define $\univi{i}$ a \deffont{universe} (at level $i$) for $i\geq 0$ as follows:
\begin{enumerate*}
\item\label{univi.i}
$\univi{i}=\hfs{i}$ for $0\leq i\leq n$, where $n$ is the choice we made in Definition~\ref{defn.n} and $\hfs{i}$ was defined in Definition~\ref{defn.hfs}. 
\item\label{univi.n+i}
$\univi{n\plus i}=\{x{\uparrow^i}\mid x\in\univi{n}\}$ for $1\leq i$.
$\uparrow$ is from Definition~\ref{defn.uparrow.and.downarrow}.
\item
We may write $\univi{}$ to refer to the totality of the model, thus $\univi{}=(\univi{i}\mid i\geq 0)$.
\end{enumerate*}
\end{defn}

\begin{rmrk}
In words we can sum up Definition~\ref{defn.univi} as follows: 
\begin{itemize*}
\item
If $0\leq i\leq n$ then an $x\in \univi{i}$ is just an element of $\hfs{i}$. 
\item
If $i>n$ then an $x\in\univi{i}$ is $u{\uparrow}^{i\minus n}$ the $(i\minus n)$-fold raising of some element $u\in\univi{n}$.
\end{itemize*} 
This is what Figure~\ref{fig.simple.diagram} illustrates.
If the reader feels it would be more principled to make the dependence on the choice of $n$ explicit then they can read `$\univi{i}$' as `$\univi{i}^n$' throughout, and no harm will come of it.
\end{rmrk}

Lemma~\ref{lemm.x.up} (especially part~\ref{x.up.down.x} of that Lemma) tells us that $\univi{n\plus i}$ is in some sense just a very simple copy of $\univi{n}$ (but `raised').\footnote{$\univi{i}$ for $i<n$ also bijects with $\univi{n}$, but using a bijection $f$ that does not preserve the native $\in$-structure (and does not preserve support).  Definition~\ref{defn.model.tin}(\ref{model.tin.n}\&\ref{model.tin.i}) treats that case.}
We make this formal:
\begin{lemm}
\label{lemm.exists.unique.u}
Suppose $i\geq 0$ and $x\in\univi{n\plus i}$.
Then:
\begin{enumerate*}
\item\label{item.exists.unique.u}
There exists a unique $u\in\univi{n}$ such that $x=u{\uparrow}^i$.
\item\label{item.arrows.bij}
$\uparrow:\univi{n\plus i}\to\univi{n\plus i\plus 1}$ and 
$\downarrow:\univi{n\plus i\plus 1}\to\univi{n\plus i}$
are bijections.
\end{enumerate*}
\end{lemm}
\begin{proof}
A fact of the construction in Definition~\ref{defn.univi}(\ref{univi.n+i}), using Lemma~\ref{lemm.x.up}.
\end{proof}

\subsection{Denotations of $y\tin x$ and $\model{\phi}$} 

Definition~\ref{defn.model.tin} is where we decide how to interpret the syntax $y\tin x$ for $y\in\univi{i}$ and $x\in\univi{i\plus 1}$.
Note that we do not set $\model{y\tin x}=(y\in x)$ --- semantically this is well-defined, but it is not what we need.
Things are much more subtle:
\begin{defn}
\label{defn.model.tin}
Suppose $y\in\univi{i}$ and $x\in\univi{i\plus 1}$.
Then we define 
$$
\model{y\tin x}\in\{\top,\bot\}
$$
depending on whether $i=n\minus 1$, $i\geq n$, or $0\leq i\leq n\minus 2$ (for $n$ the fixed but arbitrary number from Definition~\ref{defn.n}), as follows: 
\begin{enumerate*}
\item
\label{model.tin.n}
Suppose $y\in\univi{n\minus 1}$ and $x\in\univi{n}$.
We define:
$$
\model{y\tin x}= 
(f(y){\downarrow}\in x{\downarrow}{\uparrow}).
$$
For the reader's convenience we spell out the levels above:
\begin{itemize*}
\item
$y\in\univi{n\minus 1}$ so $f(y)\in\univi{n}$ and $f(y){\downarrow}\in\univi{n\minus 1}$.
\item
$x\in\univi{n}$ and $x{\downarrow}{\uparrow}\in\univi{n}$.
\end{itemize*}
Above, $f$ is the bijection between $\hfs{n\minus 1}$ and $\hfs{n}$ whose existence was noted in Corollary~\ref{corr.exists.f}.
$\uparrow$ and $\downarrow$ are an adjoint pair of (elementary) sets mappings from Definition~\ref{defn.uparrow.and.downarrow}.
\item
\label{model.tin.n+i}
Suppose $i\geq 0$ and $y\in\univi{n\plus i}$ and $x\in\univi{n\plus i\plus 1}$.
By Lemma~\ref{lemm.exists.unique.u}(\ref{item.exists.unique.u})
$y=v{\uparrow^i}$ for some unique $v\in\univi{n}$, and $x=u{\uparrow^{i\plus 1}}$ for some unique $u\in\univi{n}$.
We define:
$$
\model{y\tin x}= 
(v{\downarrow}\in u{\downarrow}{\uparrow})
.
$$
\item
\label{model.tin.i}
For $y\in\univi{i}$ and $x\in\univi{i\plus 1}$ for $0\leq i\leq n\minus 2$ we define
$$
\model{y\tin x}=\model{f(y)\in f(x)} . 
$$
\noindent Recall that $f(y)\in\univi{i\plus 1}$ and $f(x)\in\univi{i\plus 2}$, so intuitively this clause shifts $y\tin x$ up $n\minus i\minus 1$ levels until we reach $f^{n\minus i\minus 1}(y)\in\univi{n\minus 1}$ and $f^{n\minus i\minus 1}(x)\in\univi{n}$, at which point we can proceed to clause~1 above. 
\end{enumerate*}
\end{defn}

\begin{defn}[Denotation of predicates]
\label{defn.denotation}
\leavevmode
Suppose $\phi$ is a stratified predicate with $\fv(\phi)=\varnothing$ (i.e. it is closed).
It may mention constants from $\univi{}$.
We extend Definition~\ref{defn.model.tin} to a \deffont{denotation} $\model{\phi}\in\{\top,\bot\}$ as follows:
$$
\model{\tbot}=\bot
\qquad
\model{\phi\timp\phi'}=(\model{\phi}\limp\model{\phi'})
\qquad
\model{\tall a.\phi}= \compactwedge{x{\in}\univi{\level(a)}}\model{\phi[a\ssm x]}
$$
\end{defn}

\begin{rmrk}
The $\bot$, $\limp$, and $\bigwedge$ in Definition~\ref{defn.denotation} are the standard operations on the Boolean algebra $\{\top,\bot\}$ (not TST+ syntax constructors).
\end{rmrk}

\begin{xmpl}
\label{xmpl.model}
The denotation $\model{\text{-}}$ from Definitions~\ref{defn.model.tin} and~\ref{defn.denotation} is a complex design, so we give some examples.
The reader can skip this if they just want proofs, but the examples do illustrate technical points which arise in those proofs.
\begin{enumerate}
\item
Suppose $x\in\univi{n\plus 1}$ and $a\in\varsi{n}$.
We consider what it means to evaluate a predicate $\model{\tall a.a\tin x}$ asserting intuitively that $x$ is the universal set at level $n\plus 1$.
We simplify one step at a time:
\begin{itemize*}
\item
By Lemma~\ref{lemm.exists.unique.u}(\ref{item.exists.unique.u}), $x=u{\uparrow}$ for some unique $u\in\univi{n}$.
\item
Following Definitions~\ref{defn.denotation} and~\ref{defn.model.tin}(\ref{model.tin.n+i}), $\model{\tall a.a\tin x}$ when 
$\Forall{y{\in}\univi{n}}y{\downarrow}\in u{\downarrow}{\uparrow}$. 
\item
Using Lemma~\ref{lemm.arrows.adjoint} this is equivalent to 
$\Forall{y{\in}\univi{n}}y{\downarrow}{\downarrow}\in u{\downarrow}$. 
\item
By Lemma~\ref{lemm.x.up}(\ref{item.downarrow.surjective}) $\downarrow$ is surjective from $\univi{n}$ to $\univi{n\minus 1}$ to $\univi{n\minus 2}$, so this simplifies to 
$\Forall{y'{\in}\univi{n\minus 2}}y'\in u{\downarrow}$.
\item 
So $\model{\tall a.a\tin u{\uparrow}}$ holds when $u{\downarrow}=\univi{n\minus 1}$.
This is also equivalent to $x\downone{2}=\univi{n\minus 1}$.
\end{itemize*}
The reader might have expected $\model{\tall a.a\tin u{\uparrow}}$ to mean $u{\uparrow}=\univi{n\plus 1}$ rather than $u{\downarrow}=\univi{n\minus 1}$.
This is a feature, not a bug: but for now let us just note this fact and continue to more examples.
\item
Suppose $x\in\univi{n}$ and $a\in\varsi{n\minus 1}$.
Again, we evaluate $\model{\tall a.a\tin x}$.
\begin{itemize*}
\item
Following Definitions~\ref{defn.denotation} and~\ref{defn.model.tin}(\ref{model.tin.n}), 
$\model{\tall a.a\tin x}$ when 
$\Forall{y{\in}\univi{n\minus 1}}f(y){\downarrow}\in x{\downarrow}{\uparrow}$. 
\item
Using Lemma~\ref{lemm.arrows.adjoint} this is equivalent to
$\Forall{y{\in}\univi{n\minus 1}}f(y){\downarrow}{\downarrow}\in x{\downarrow}$. 
\item
By Corollary~\ref{corr.exists.f} $f$ bijects $\univi{n\minus 1}$ with $\univi{n}$, so we can simplify to 
$\Forall{y'{\in}\univi{n}}y'{\downarrow}{\downarrow}\in x{\downarrow}$. 
\item
By Lemma~\ref{lemm.x.up}(\ref{item.downarrow.surjective}) $\downarrow$ is surjective, so this simplifies to 
$\Forall{y''{\in}\univi{n\minus 2}}y''\in x{\downarrow}$.
\item 
So $\model{\tall a.a\tin x}$ holds when $x{\downarrow}=\univi{n\minus 1}$.
\end{itemize*}
\item
We leave it to the reader to check that:
\begin{enumerate*}
\item
When $x\in\univi{n\plus i\plus 1}$ and $a\in\varsi{n\plus i}$ for $i\geq 0$ then $\model{\tall a.a\tin x}$ is equivalent to $x\downone{i\plus 2}=\univi{n\minus 1}$.
\item
When $x\in\univi{i\plus 1}$ and $a\in\varsi{i}$ for $0\leq i\leq n\minus 1$ then $\model{\tall a.a\tin x}$ is equivalent to $f^{n\minus i\minus 1}(x){\downarrow}=\univi{n\minus 1}$.
\end{enumerate*}
\item
What to make of the examples so far?

Intuitively, we can think of $\model{\tall y.y\tin x}$ as asserting that $x{\downarrow}$ is equal to $\univi{n\minus 1}$, for $x$ \emph{as it presents itself at $\univi{n}$}.
If $x$ is lower than $n$ then we use $f$ to shift it up; if $x$ is higher than $n$ then we use $\downarrow$ to shift it down.

\item
Consider $x,x'\in\univi{n}$ and $a\in\varsi{n\minus 1}$.
We unpack what $\model{\tall a.a\tin x\tiff a\tin x'}$ means.

We follow Definitions~\ref{defn.denotation} and~\ref{defn.model.tin}(\ref{model.tin.n}) and simplify to
$$
\Forall{y{\in}\univi{n\minus 1}}f(y){\downarrow}\in x{\downarrow}{\uparrow} \liff f(y){\downarrow}\in x'{\downarrow}{\uparrow} .
$$
After simplifications using Lemmas~\ref{lemm.arrows.adjoint} and~\ref{lemm.x.up}(\ref{item.downarrow.surjective}) and Corollary~\ref{corr.exists.f}, as we already saw in the examples above, this reduces to
$$
x{\downarrow}= x'{\downarrow} .
$$
We will see two generalisations of this property below: 
\begin{itemize*}
\item
there is (most obviously) Proposition~\ref{prop.extensionality} (extensionality), whose proof visibly mirrors this example; but also
\item
there is Lemma~\ref{lemm.updown}, which shows that $x$ and $x{\downarrow}{\uparrow}$ are \emph{always} extensionally equal with respect to our denotation $\model{\text{-}}$.
This property is a feature, and is needed for Proposition~\ref{prop.comprehension}.\footnote{It does \emph{not} follow that $x$ is extensionally equivalent to $x\downone{i}{\uparrow}^i$ for arbitrary $i$.  This would be wrong.  What does hold is $(x{\downarrow}{\uparrow}){\downarrow}{\uparrow} = x{\downarrow}{\uparrow}$ using Lemma~\ref{lemm.x.up}(\ref{x.up.down.x}).  To use some jargon: $x\mapsto x{\downarrow}{\uparrow}$ is idempotent.} 
\end{itemize*}
\end{enumerate}
A theme of the examples above is this: the $\in$-structure of the denotation is determined by \emph{`$\univi{n}$-shifted-down-once-with-$\downarrow$'} 
(the reader who wants to see this intuition made precise and brought to bear in full technical detail, can look to the proof of Lemma~\ref{lemm.updown}).
If an element $x$ is lower than $n$ then we use $f$ to (by Corollary~\ref{corr.exists.f}) bijectively shift $x$ up to $n$, and then we apply $\downarrow$; if $x$ is higher than $n$ then we use $\downarrow$ to (by Lemma~\ref{lemm.exists.unique.u}(\ref{item.arrows.bij})) bijectively shift $x$ down to $n$, and we apply $\downarrow$. 
\end{xmpl}

\subsection{Shift and typical ambiguity}

\begin{defn}
\label{defn.shift}
Suppose $x\in\univi{i}$ for $i\geq 0$.
Then define $x^\plus\in\univi{i\plus 1}$, read \deffont{shift} $x$, as follows:
\begin{enumerate*}
\item\label{item.shift.i}
$x^\plus = f(x)$ if $x\in\univi{i}$ for $0\leq i\leq n\minus 1$.
\item\label{item.shift.n+i}
$x^\plus = x{\uparrow}$ if $x\in\univi{n\plus i}$ for $0\leq i$.
\end{enumerate*}
In words: at levels $i<n$ we shift up using the bijection $f$ from Corollary~\ref{corr.exists.f}; at levels $i\geq n$ we shift up using the bijection that is built in to the model at higher levels, as per Lemma~\ref{lemm.exists.unique.u}.
\end{defn}

\begin{lemm}
\label{lemm.shift.bij}
The shift mapping $x\mapsto x^\plus$ bijects $\univi{i}$ with $\univi{i\plus 1}$, for every $i\geq 0$.
\end{lemm}
\begin{proof}
$f:\univi{i}\to\univi{i\plus 1}$ for $0\leq i\leq n\minus 1$ is a bijection by Corollary~\ref{corr.exists.f}.
$\uparrow:\univi{n\plus 1}\to\univi{n\plus i\plus 1}$ for $0\leq i$ is a bijection by Lemma~\ref{lemm.exists.unique.u}(\ref{item.arrows.bij}).
\end{proof}

We now reprise Notation~\ref{nttn.plus}(\ref{item.upsome}), but at a slightly higher level of sophistication because we have constructed our set of constants:
\begin{defn}
\label{defn.shift.phi}
Make a fixed but arbitrary choice of bijection $a\mapsto a^\plus$ between the variable symbols at levels $i$ and $i\plus 1$, for each $i\geq 0$.
We use this bijection to extend the shift action from $\univi{}$ to closed predicates (possibly mentioning constants) by defining $\phi^\plus$ to be a copy of $\phi$ in which 
\begin{itemize*}
\item
each constant $x$ appearing in $\phi$ is replaced by $x^\plus$, and 
\item
each variable symbol $a$ with $\level(a)=i$ appearing in $\phi$ is replaced 
with our fixed but arbitrary bijective choice of variable symbol $a^\plus$ such that $\level(a^\plus)=i\plus 1$. 
\end{itemize*}
For example: if $x\in\univi{n\plus 1}$ and $\level(a)=\level(b)=n$ then $x^\plus=x{\uparrow}$ and $(\tall a.\texi b.b\tin a\timp \tneg (b\tin x))^\plus$ is $\tall a^\plus.\texi b^\plus.b^\plus\tin a^\plus\timp \tneg (b^\plus\tin x{\uparrow})$. 
An inductive definition is not hard to write out.
\end{defn}

\begin{rmrk}
\label{rmrk.TA.consts}
We make two elementary observations about Proposition~\ref{prop.TA} (typical ambiguity) before going into the details:
\begin{enumerate*}
\item
The typical ambiguity axiom \rulefont{TA_\phi} in Figure~\ref{fig.TZT.axioms} requires us to consider closed $\phi$, without constants.
Proposition~\ref{prop.TA} will prove something stronger: that $\model{\phi^\plus}=\model{\phi}$ for any closed $\phi$ (possibly mentioning constants from $\univi{}$).
Clearly, if typical ambiguity possibly with constants is valid, then typical ambiguity without constants follows as a special case.

How does this extra generality help the proofs?
It helps because in this paper we care about denotation in \emph{one specific model} $\univi{}$.
$\univi{}$ has a shift symmetry built in (using $f$ and $\uparrow$; see Definition~\ref{defn.shift}), and furthermore, the denotation in Definition~\ref{defn.denotation} introduces constant symbols from $\univi{}$.
THus, even if we start from a closed predicate without constants, an induction over syntax is likely to lead us to consider closed predicates \emph{with} constants. 
Thus for our specific constructionns, typical ambiguity on predicates with constants is both meaningful and necessary.
\item
An intuition for why Proposition~\ref{prop.TA} is plausible, is that Definition~\ref{defn.model.tin} interprets $y\tin x$ as $x$ and $y$ presents themselves at $\univi{n}$ (see the discussion in Subsection~\ref{subsect.some.intuitions}).
If $y$ and $x$ are not in $\univi{n}$, then they get shifted so they are (up using $f$; down using $\downarrow$).
In accordance with this intuition, $\model{\phi^+}=\model{\phi}$ holds because everything in $\phi^+$ and $\phi$ gets shifted to $\univi{n}$ anyway.
\end{enumerate*}
\end{rmrk}

\begin{prop}[\bf Typical ambiguity]
\label{prop.TA}
Suppose $\phi$ is a closed predicate (possibly mentioning constants).
Then 
$$
\model{\phi^\plus}=\model{\phi}.
$$
\end{prop}
\begin{proof}
We reason by induction on $\phi$:
\begin{itemize}
\item
\emph{The case of $\tall a.\phi$ where $\level(a)=i$.}\quad
We reason as follows:
$$
\begin{array}{r@{\ }l@{\quad}l}
\ment \tall a^\plus.\phi^\plus
\liff&
\Forall{x'{\in}\univi{i\plus 1}}\ment\phi^\plus[a^\plus\ssm x']
&\text{Definition~\ref{defn.denotation}}
\\
\liff&
\Forall{x{\in}\univi{i}}\ment \phi^\plus[a^\plus\ssm x^\plus]
&\text{Lemma~\ref{lemm.shift.bij}}
\\
\liff&
\Forall{x{\in}\univi{i}}\ment \phi[a\ssm x]
&\text{Ind. hyp. on $\phi[a\ssm x]$}
\\
\liff&
\ment\tall a.\phi
&\text{Definition~\ref{defn.denotation}}
\end{array}
$$
\item
\emph{The cases of $\tbot$ and $\phi\timp\phi'$.}\quad
By routine logical manipulations.
\item
\emph{The case of $y\tin x$ where $y\in\univi{i}$ and $x\in\univi{i\plus 1}$ for $0\leq i\leq n\minus 1$.}\quad

By Definition~\ref{defn.shift} $y^\plus=f(y)\in\univi{i\plus 1}$ and $x^\plus=f(x)\in\univi{i\plus 2}$, and combining with Definition~\ref{defn.shift.phi} we have that
$$
\model{(y\tin x)^\plus}=\model{f(y)\tin f(x)}.
$$
Meanwhile, by Definition~\ref{defn.model.tin}(\ref{model.tin.i}) also $\model{y\tin x}=\model{f(y)\tin f(x)}$, so we are done.
\item
\emph{The case of $y\tin x$ where $y\in\univi{n\minus 1}$ and $x\in\univi{n}$.}\quad

By Definition~\ref{defn.shift} $y^\plus=f(y)\in\univi{n}$ and $x^\plus=x{\uparrow}\in\univi{n\plus 1}$.
So we must check that 
$$
\model{y\tin x}=\model{f(y)\tin x{\uparrow}}. 
$$
We unpack clauses~\ref{model.tin.n} and~\ref{model.tin.n+i} of Definition~\ref{defn.model.tin} and see that
$$
\model{y\tin x}= (f(y){\downarrow}\in x{\downarrow}{\uparrow})
\quad\text{and}\quad
\model{f(y)\tin x{\uparrow}}= (f(y){\downarrow} \in x{\downarrow}{\uparrow})  .
$$
These conditions are identical, so we are done.
\item
\emph{The case of $y\tin x$ where $y\in\univi{n\plus i}$ and $x\in\univi{n\plus i\plus 1}$ for $i\geq 0$.}\quad

Unpacking Definition~\ref{defn.shift.phi}, $(y\tin x)^\plus = (y{\uparrow}\tin x{\uparrow})$.

By Lemma~\ref{lemm.exists.unique.u}(\ref{item.exists.unique.u})
$y=v{\uparrow}^i$ for some unique $v\in\univi{n}$, and $x=u{\uparrow}^{i\plus 1}$ for some unique $u\in\univi{n}$, and by Definition~\ref{defn.model.tin}(\ref{model.tin.n+i})
$$
\model{y\tin x} = (v{\downarrow}\in u{\downarrow}{\uparrow}) = \model{y{\uparrow}\tin x{\uparrow}} , 
$$
so we are done.
\qedhere\end{itemize} 
\end{proof}

\subsection{Extensionality}

\begin{prop}[\bf Extensionality]
\label{prop.extensionality}
Let $a$, $a'$, $b$, and $c$ be variables such that $\level(a)=\level(a')=i$ and $\level(b)=i\minus 1$ and $\level(c)=i\plus 1$.
Then
$$
\model{\tall a,a'.(\tall b.b\tin a\tiff b\tin a')\timp \tall c.a\tin c\tiff a'\tin c}=\top .
$$ 
\end{prop}
\begin{proof}
Suppose without loss of generality using Proposition~\ref{prop.TA} that $i=n\plus 1$; 
so $\level(b)=n$ and $\level(a)=\level(a')=n\plus 1$ and $\level(c)=n\plus 2$.

Choose any $x,x'\in\univi{n\plus 1}$.
Following Definition~\ref{defn.denotation} it suffices to proceed as follows: 
\begin{itemize*}
\item
Instantiate $a$ to $x\in\univi{n\plus 1}$ and $a'$ to $x'\in\univi{n'\plus 1}$.
\item
Suppose $\model{\tall b.b\tin x\tiff b\tin x'}=\top$.
\item
Seek to prove $\model{\tall c.x\tin c\tiff x'\tin c}=\top$.
\end{itemize*}
So fix $x$ and $x'$ and suppose $\model{\tall b.b\tin x\tiff b\tin x'}=\top$.
Unpacking Definition~\ref{defn.denotation}, we have 
$$
\Forall{y\in\univi{n}}\model{y\tin x}\liff \model{y\tin x'} .
$$ 
By Lemma~\ref{lemm.exists.unique.u}(\ref{item.exists.unique.u}) (since $x,x'\in\univi{n\plus 1}$) 
$$
x=u{\uparrow}
\quad\text{and}\quad
x'=u'{\uparrow}
$$ 
for some unique $u,u'\in\univi{n}$.
From Definition~\ref{defn.model.tin}(\ref{model.tin.n+i}) we conclude that
$$
\Forall{y\in\univi{n}}y{\downarrow}\in u{\downarrow}{\uparrow}\liff y{\downarrow}\in u'{\downarrow}{\uparrow} .
$$
Now by Lemma~\ref{lemm.x.up}(\ref{item.downarrow.surjective}) $\downarrow$ is surjective from $\univi{n}$ to $\univi{n\minus 1}$, so 
$$
\Forall{v\in\univi{n\minus 1}}v\in u{\downarrow}{\uparrow}\liff v\in u'{\downarrow}{\uparrow} 
\quad\text{and by sets extensionality}\quad
u{\downarrow}{\uparrow}=u'{\downarrow}{\uparrow}.
$$
By Lemma~\ref{lemm.x.up}(\ref{item.uparrow.injective}) $\uparrow$ is injective and we conclude that
$$
u{\downarrow}= u'{\downarrow} .
$$
Now consider $z\in\univi{n\plus 2}$.
By Lemma~\ref{lemm.exists.unique.u}(\ref{item.exists.unique.u}) $z=w{\uparrow^2}$ for some unique $w\in\univi{n}$.
By Definition~\ref{defn.model.tin}(\ref{model.tin.n+i}) it just remains to check that 
$$
u{\downarrow}\in w{\downarrow}{\uparrow} \liff u'{\downarrow}\in w{\downarrow}{\uparrow} .
$$
But this is immediate since $u{\downarrow}=u'{\downarrow}$, so we are done.
\end{proof}

\subsection{The permutation action on predicates}

\begin{rmrk}
It will be useful to define a permutation action on predicates.

It might seem easy to do this by extending the action on $\univi{}$ and letting $\pi\act\phi$ be a copy of $\phi$ in which each constant $x$ appearing in $\phi$ is replaced by $\pi\act x$. 
But this is not what we need, because of how the bijection $f$ is used to define the $\in$-structure of our model as defined (for example) in Definition~\ref{defn.model.tin}(\ref{model.tin.n}).
This is also where much of the work in proof of Proposition~\ref{prop.TA} happens.

So the pointwise permutation action $\pi\act x$ is not quite what we need on $\univi{i}$ for $0\leq i\leq n\minus 1$.
We now define an `adjusted' permutation action that takes our use of $f$ into account, as follows: 
\end{rmrk}

\begin{defn}
\label{defn.permutation.action}
Suppose $x\in\univi{i}$ for $i\geq 0$ and $\pi\in\tf{Perm}$ is a permutation.
Then define $\pi\mact x$ the \deffont{adjusted permutation action} by: 
\begin{enumerate*}
\item\label{item.mact.n+i}
$\pi\mact x = \pi\act x$ if $x\in\univi{n\plus i}$ for $i\geq 0$.
\item\label{item.mact.i}
$\pi\mact x = f^\mone(\pi\mact (f(x)))$ if $x\in\univi{i}$ for $0\leq i\leq n\minus 1$.
\end{enumerate*} 
Extend the action $\mact$ to predicates $\phi$ (possibly with constants; possibly with free variables) such that $\pi\mact\phi$ is that predicate obtained by replacing each constant $x$ appearing in $\phi$, with $\pi\mact x$.
\end{defn} 

\begin{rmrk}
Clause~\ref{item.mact.i} in Definition~\ref{defn.permutation.action} defines $\pi\mact x$ in terms of $\pi\mact f(x)$, for $x\in\univi{i}$.
This is inductive on $n\minus i$; when this quantity reduces to $0$ we exit clause~\ref{item.mact.i} and move to clause~\ref{item.mact.n+i}.
It is not hard to see that an equivalent formulation of clause~\ref{item.mact.i} is that 
$$
\pi\mact x = f^{i\minus n}(\pi\act (f^{n\minus i}(x)))
$$ 
for $x\in\univi{i}$ for $0\leq i\leq n\minus 1$.

Note also that if $x\in\univi{n\plus i}$ then $x=u{\uparrow^i}$ for some $u\in\univi{n}$ by Lemma~\ref{lemm.exists.unique.u}(\ref{item.exists.unique.u}).
Also by Lemma~\ref{lemm.arrows.equivariant} $\pi$ commutes with $\uparrow$ and $\downarrow$, so that $\pi\act(u{\uparrow^i})=(\pi\act u){\uparrow^i}$. 

So we can think of $\pi\mact x$ in Definition~\ref{defn.permutation.action} as being all about $\pi$ acting on $x$ \emph{as $x$ presents itself at} $\univi{n}$.
If $x$ is lower than $n$ then we use $f$ to raise it, act, then use $f^\mone$ to go back; if $x$ is higher than $n$ then we use $\downarrow$ to lower it, act, then use $\uparrow$ to go back.
\end{rmrk}

A test that we have got our definition of $\mact$ right is that the permutation action should commute with shift:
\begin{lemm}
\label{lemm.pi.mact.i}
Suppose $x\in\univi{i}$ for $0\leq i\leq n\minus 1$ and $\pi\in\tf{Perm}$ is a permutation.
Then
$$
(\pi\mact x)^\plus = \pi\mact f(x) = \pi\mact (x^\plus) . 
$$ 
\end{lemm}
\begin{proof}
The right-hand equality is immediate because $x^\plus=f(x)$ by Definition~\ref{defn.shift}(1).

For the left-hand equality we reason as follows:
$$
\begin{array}{r@{\ }l@{\qquad}l}
(\pi\mact x)^\plus 
=&
(f^\mone(\pi\mact f(x)))^\plus
&\text{Definition~\ref{defn.permutation.action}(\ref{item.mact.i})}
\\
=&
f(f^\mone(\pi\mact f(x)))
&\text{Definition~\ref{defn.shift}(\ref{item.shift.i})}
\\
=&
\pi\mact f(x) 
&f\circ f^\mone=\f{id}
\end{array}
$$
\end{proof}

\begin{corr}
\label{corr.pi.mact.plus}
Suppose $x\in\univi{i}$ for $i\geq 0$ and suppose $\pi\in\tf{Perm}$ is a permutation. 
Then
$$
\pi\mact(x^\plus) = (\pi\mact x)^\plus .
$$
\end{corr}
\begin{proof}
If $0\leq i\leq n\minus 1$ then this is immediate from Lemma~\ref{lemm.pi.mact.i}.

If $i\geq n$ then we reason as follows:
$$
\begin{array}{r@{\ }l@{\qquad}l}
\pi\mact (x^\plus) 
=& \pi\mact (x{\uparrow})
&\text{Definition~\ref{defn.shift}(\ref{item.shift.n+i})}
\\
=& \pi\act (x{\uparrow})
&\text{Definition~\ref{defn.permutation.action}(\ref{item.mact.n+i})}
\\
=& (\pi\act x){\uparrow}
&\text{Lemma~\ref{lemm.arrows.equivariant}}
\\
=& (\pi\mact x){\uparrow}
&\text{Definition~\ref{defn.permutation.action}(\ref{item.mact.n+i})}
\\
=& (\pi\mact x)^\plus
&\text{Definition~\ref{defn.shift}(\ref{item.shift.n+i})}
\end{array}
$$
\end{proof}

\begin{rmrk}
Note for experts who are familiar with nominal terms: $\mact$ is a bit simplistic in that it does not
suspend permutations on free variables (i.e. $\pi\mact a$ is not syntax, for $a$ a variable symbol).
To add this is very feasible --- it is what the syntax of \emph{nominal terms} does with its \emph{suspended permutations}, as used in \emph{nominal unification} and \emph{rewriting}~\cite{gabbay:nomu-jv,gabbay:nomr-jv} --- but we do not need suspended permutation syntax here. 
An expert in nominal terms might want some explanation of why this is so.

Here, we are evaluating validity of predicates with respect to one particular (and very special) model $\univi{}$.
We know what values variables range over and can include them as constants in syntax.
This mostly avoids dealing with open predicates and, now, suspended permutations.
\end{rmrk}

\begin{prop}
\label{prop.shift.equivar}
Suppose $\phi$ is a closed predicate (it may mention constants) and suppose $\pi\in\tf{Perm}$ is a permutation.
Then
$$
\pi\mact(\phi^\plus) = (\pi\mact\phi)^\plus .
$$
\end{prop}
\begin{proof}
By induction on $\phi$.\footnote{Actually, not quite, since a subterm of a closed predicate is not a closed predicate.  We can deal with this in three ways: strengthen the induction to account for it; or use a fudge encoding and instantiate each quantified variable to a freshly-chosen constant as we induct into the term --- or, we can note that what we actually need for the proofs is $\model{\pi\mact(\phi^\plus)}=\model{(\pi\mact\phi)^\plus}$, and \emph{this} we can prove by inducting just over closed predicates, since quantifiers instantiate.} 
The cases of $\tbot$, $\timp$, and $\tall$ are routine.
The case of $y\tin x$ is just from Corollary~\ref{corr.pi.mact.plus}.
\end{proof}

\subsection{Comprehension}
\label{subsect.comprehension}

We are aiming for Proposition~\ref{prop.comprehension} (validity of comprehension), but we need to build some technical machinery first.

\begin{lemm}
\label{lemm.updown}
Suppose that:
\begin{itemize*}
\item
$\phi$ is a predicate (possibly mentioning constants), and 
$\fv(\phi)=\{a\}$. 
\item
$\level(a)=n$ --- note that $a$ has level $n$ specifically for $n$ as chosen in Definition~\ref{defn.n}; the proof will depend on this.
\item
$x\in\univi{n}$.
\end{itemize*}
Then 
$$
\model{\phi[a\ssm x]} =
\model{\phi[a\ssm x{\downarrow}{\uparrow}]} .
$$
\end{lemm}
\begin{proof}
We reason by induction on $\phi$.
The cases of $\tbot$, $\timp$, and $\tall$ are routine from Definition~\ref{defn.denotation}.
We sketch the case for $\tall$ to show that this is indeed routine:
\begin{multline*}
\model{(\tall b.\phi)[a\ssm x]} 
\stackrel{\text{fact}}=
\model{\tall b.(\phi[a\ssm x])}
\stackrel{D\ref{defn.denotation}}=
\compactwedge{y{\in}\univi{\level(b)}}\model{\phi[a\ssm x][b\ssm y]}
\stackrel{\text{fact}}=
\\
\compactwedge{y{\in}\univi{\level(b)}}\model{\phi[b\ssm y][a\ssm x]}
\stackrel{ind. hyp.}=
\compactwedge{y{\in}\univi{\level(b)}}\model{\phi[b\ssm y][a\ssm x{\downarrow}{\uparrow}]}
\stackrel{\text{same reasoning in reverse}}{=\dots=}
\model{(\tall b.\phi)[a\ssm x{\downarrow}{\uparrow}]} 
\end{multline*}
The interesting case is the case for $\tin$.
Due to restrictions on levels, we have three cases: 
\begin{enumerate*}
\item
$\phi=(a\tin y)$ for $y\in\univi{n\plus 1}$. 
\item
$\phi=(y\tin a)$ for $y\in\univi{n\minus 1}$. 
\item
$\phi=(y\tin y')$ for $y\in\univi{i}$ and $y'\in\univi{i\plus 1}$ for $i\geq 0$. 
\end{enumerate*}
We consider each in turn:
\begin{enumerate}
\item
\emph{The case of $\phi=(a\tin y)$ for $x\in\univi{n}$ and $y\in\univi{n\plus 1}$.}\quad

By Lemma~\ref{lemm.exists.unique.u}(\ref{item.exists.unique.u})
 $y=v{\uparrow}$ for some unique $v\in\univi{n}$, and by Definition~\ref{defn.model.tin}(\ref{model.tin.n+i}) 
$$
\model{x\tin y}=(x{\downarrow}\in v{\downarrow}{\uparrow}). 
$$
By Lemma~\ref{lemm.x.up}(\ref{x.up.down.x}) $x{\downarrow}=(x{\downarrow}){\uparrow}{\downarrow}=(x{\downarrow}{\uparrow}){\downarrow}$, and the result follows.
\item
\emph{The case of $\phi=(y\tin a)$ for $y\in\univi{n\minus 1}$ and $x\in\univi{n}$.}\quad

By Definition~\ref{defn.model.tin}(\ref{model.tin.n}) 
$$
\model{y\tin x}=(f(y){\downarrow}\in x{\downarrow}{\uparrow}).
$$
Using Lemma~\ref{lemm.x.up}(\ref{x.up.down.x}) 
$x{\downarrow}{\uparrow}=((x{\downarrow}){\uparrow}{\downarrow}){\uparrow}= (x{\downarrow}{\uparrow}){\downarrow}{\uparrow}$, and the result follows.
\item
\emph{The case of $\phi=(y\tin y')$ for $y\in\univi{i}$ and $y'\in\univi{i\plus 1}$ for $i\geq 0$.}\quad

There is nothing to prove here, since $(y\tin y')[a\ssm x]=(y\tin y')=(y\tin y')[a\ssm x{\downarrow}{\uparrow}]$. 
\qedhere\end{enumerate}
\end{proof}

\begin{prop}[Equivariance]
\label{prop.equivariance}
Suppose $\phi$ is a closed predicate (without free variables; possibly with constants), and suppose $\pi\in\tf{Perm}$ is a permutation.
Then
$$
\model{\phi} = \model{\pi\mact\phi} .
$$
\end{prop}
\begin{proof}
First, note from Propositions~\ref{prop.TA} and~\ref{prop.shift.equivar} that
$$
\model{\phi}=\model{\phi^\plus}
\quad\text{and}\quad
\model{\pi\mact\phi}=\model{(\pi\mact\phi)^\plus}=\model{\pi\mact(\phi^\plus)} .
$$
Thus we can shift levels up and assume without loss of generality that every constant and variable symbol in $\phi$ has level at least $n$.
This eliminates some cases in the proof below. 

We now reason by induction on $\phi$.
We consider each case in turn:
\begin{itemize}
\item
\emph{The case of $\tall a.\phi$.}\quad

We reason as follows:
$$
\begin{array}{r@{\ }l@{\qquad}l}
\model{\pi\mact\tall a.\phi}
=&
\model{\tall a.\pi\mact\phi}
&\text{Definition~\ref{defn.permutation.action}}
\\
=&
\bigwedge_{x{\in}\univi{\level(a)}}\model{(\pi\mact\phi)[a\ssm x]}
&\text{Definition~\ref{defn.denotation}}
\\
=&
\bigwedge_{x{\in}\univi{\level(a)}}\model{(\pi\mact\phi)[a\ssm \pi\mact x]}
&\text{Lem.~\ref{lemm.pi.pointwise}(\ref{item.pi.bij}), Def.~\ref{defn.permutation.action}(\ref{item.mact.n+i}), $\level(a)\geq n$}
\\
=&
\bigwedge_{x{\in}\univi{\level(a)}}\model{\pi\mact(\phi[a\ssm x])}
&\text{Fact, $\level(a)\geq n$}
\\
=&
\bigwedge_{x{\in}\univi{\level(a)}}\model{\phi[a\ssm x]}
&\text{Ind. hyp.}
\\
=&
\model{\tall a.\phi}
&\text{Definition~\ref{defn.denotation}}
\end{array}
$$
\item
\emph{The cases of $\tbot$ and $\timp$}\quad \dots are routine.
\item
\emph{The case of $y\tin x$ where $y\in\univi{n\minus 1}$ and $x\in\univi{n}$.}\quad

We do not need to consider this because we assumed $\level(y)\geq n$.
\item
\emph{The case of $y\tin x$ where $y\in\univi{n\minus i}$ and $x\in\univi{n\minus i\plus 1}$ for $2\leq i<n$.}\quad

We do not need to consider this because we assumed $\level(y)\geq n$.
\item
\emph{The case of $y\tin x$ where $y\in\univi{n\plus i}$ and $x\in\univi{n\plus i\plus 1}$ for $0\leq i$.}\quad

By Lemma~\ref{lemm.exists.unique.u}(\ref{item.exists.unique.u}) $y=v{\uparrow}^i$ for some unique $v\in\univi{n}$, and $x=u{\uparrow}^{i\plus 1}$ for some unique $u\in\univi{n}$.
We reason as follows:
$$
\begin{array}{r@{\ }l@{\qquad}l}
\model{\pi\mact (y\tin x)} 
=& 
\model{(\pi\act y)\tin(\pi\act x)} 
&\text{Definition~\ref{defn.permutation.action}, $\level(y)\geq n$}
\\
=& 
\model{(\pi\act (v{\uparrow}^i)\tin(\pi\act (u{\uparrow}^{i\plus 1}))} 
&\text{Fact}
\\
=& 
\model{((\pi\act v){\uparrow}^i)\tin((\pi\act u){\uparrow}^{i\plus 1})} 
&\text{Lemma~\ref{lemm.arrows.equivariant}}
\\
=& 
(\pi\act v){\downarrow}\in (\pi\act u){\downarrow}{\uparrow} 
&\text{Definition~\ref{defn.model.tin}(\ref{model.tin.n+i})}
\\
=& 
\pi\act (v{\downarrow})\in \pi\act (u{\downarrow}{\uparrow}) 
&\text{Lemma~\ref{lemm.arrows.equivariant}}
\\
=&
v{\downarrow}\in u{\downarrow}{\uparrow}
&\text{Lemma~\ref{lemm.pi.pointwise}(\ref{item.pi.bij})}
\\
=&
\model{y\tin x} 
&\text{Definition~\ref{defn.model.tin}(\ref{model.tin.n+i})}
\end{array}
$$
\qedhere\end{itemize} 
\end{proof}

The adjusted permutation action $\pi\mact x$ from Definition~\ref{defn.permutation.action} has a corresponding adjusted notion of support $\adjsupp(x)$:
\begin{defn}
\label{defn.modified.supp}
Suppose $x\in\univi{i}$ for $0\leq i$.
Define $\adjsupp(x)$ the \deffont{adjusted support} of $x$ as follows:
\begin{enumerate*}
\item
If $0\leq i\leq n\minus 1$ then $\adjsupp(x)=\adjsupp(f(x))$.
\item
If $i\geq n$ then $\adjsupp(x)=\supp(x)$.
\end{enumerate*}
If $\phi$ is a closed predicate (it may mention constants) then we define
$$
\adjsupp(\phi)=\bigcup\{\adjsupp(x) \mid x\in\consts(\phi)\}.
$$
($\consts(\phi)$ denotes the set of constants appearing in $\phi$; see Notation~\ref{nttn.plus}(\ref{item.plus.closed}).)
\end{defn}

\begin{rmrk}
Clause~1 of Definition~\ref{defn.modified.supp} is well-defined by an induction on $n\minus i$; when this quantity reaches $0$ we exit to clause~2.
It is easy to see that:
\begin{itemize*}
\item
If $0\leq i\leq n\minus 1$ then $\adjsupp(x)=\supp(f^{n\minus i}(x))$.
\item
If $i\geq n$ then $\adjsupp(x)=\supp(x)$.
\end{itemize*}
Thus, the adjusted support $\adjsupp(x)$ works on $x\in\univi{i}$ for $0\leq i\leq n\minus 1$ by shifting $x$ up using $f$ until it is in $\univi{n}$, and evaluating $\supp(x)$ at $\univi{n}$.
Also, it follows from Lemmas~\ref{lemm.exists.unique.u}(\ref{item.arrows.bij}) and~\ref{lemm.arrows.equivariant} that if $x\in\univi{n\plus i}$ for $i\geq 1$ then $\supp(x)=\supp(x{\downarrow})$.
So intuitively we can think of $\adjsupp(x)$ as being `the support of $x$ as it presents itself in $\univi{n}$'.
\end{rmrk}

\begin{lemm}
\label{lemm.fix.supp'.x}
Suppose $x\in\univi{i}$ for $i\geq 0$ and $\pi\in\tf{Perm}$ is a permutation.
Then if $\pi\in\fix(\adjsupp(x))$ (Definition~\ref{defn.fix.S}) then $\pi\mact x=x$.
\end{lemm}
\begin{proof}
Suppose $\pi\in\fix(\adjsupp(x))$.
We consider cases:
\begin{enumerate*}
\item
If $i\geq n$ then $\adjsupp(x)=\supp(x)$ and $\pi\mact x=\pi\act x$ and we just apply Lemma~\ref{lemm.pi.pointwise}(\ref{item.supp.supports}).
\item
If $0\leq i\leq n\minus 1$ then $\adjsupp(x)=\adjsupp(f(x))$ and $\pi\mact x = f^\mone(\pi\mact f(x))$; so it would suffice to prove that $\pi\mact f(x)=f(x)$.
Now $f(x)\in\univi{i\plus 1}$ and the result follows by a routine induction on $n\minus i$, using as base case the previous clause with $i=n$ (i.e. $n\minus i=0$).
\qedhere\end{enumerate*}
\end{proof}

\begin{corr}
\label{corr.sub.pi.mact.x}
Suppose that:
\begin{itemize*}
\item
$\phi$ is a predicate (it may mention constants) and $\fv(\phi)=\{a\}$. 
\item
$x\in\univi{\level(a)}$.
\item
$\pi\in\tf{Perm}$ is a permutation. 
\end{itemize*}
Then
$$
\pi\in\fix(\adjsupp(\phi))
\quad\text{implies}\quad
\model{\phi[a\ssm x]}=\model{\phi[a\ssm \pi\mact x]}.
$$
\end{corr}
\begin{proof}
By Proposition~\ref{prop.equivariance} $\model{\phi[a\ssm x]}=\model{\pi\mact(\phi[a\ssm x])}$.
By Lemma~\ref{lemm.fix.supp'.x} $\pi\mact y=y$ for every constant $y$ appearing in $\phi$, so that $\pi\mact(\phi[a\ssm x])=\phi[a\ssm \pi\mact x]$, and thus $\model{\pi\mact(\phi[a\ssm x])}=\model{\phi[a\ssm \pi\mact x]}$ as required. 
\end{proof}

\begin{corr}
\label{corr.X.in.Un}
Suppose $\phi$ is a predicate (it may mention constants) and $\fv(\phi)=\{a\}$ and $\level(a)=n$.
Then 
$$
X=\{x{\downarrow} \mid x\in\univi{n},\ \model{\phi[a\ssm x{\downarrow}{\uparrow}]}=\top \}\subseteq\univi{n\minus 1}
$$
is supported (in the sense of Definition~\ref{defn.S.supports.x}(\ref{item.S.supports.x})) by $\adjsupp(\phi)$, and thus (since by Definition~\ref{defn.univi}(\ref{univi.i}) $\univi{n}$ is the set of finitely supported subsets of $\univi{n\minus 1}$) we have that 
$$
X\in\univi{n} .
$$
\end{corr}
\begin{proof}
We must show that $\pi\act X=X$ for every $\pi\in\fix(\adjsupp(\phi))$.

So suppose $\pi\in\fix(\adjsupp(\phi))$; note that it follows that also $\pi^\mone\in\fix(\adjsupp(\phi))$.
Then:
$$
\begin{array}{r@{\ }l@{\qquad}l}
\pi\act X
=&
\{\pi\act (x{\downarrow}) \mid x\in\univi{n},\ \model{\phi[a\ssm x{\downarrow}{\uparrow}]}=\top\}
&\text{Pointwise action (Definition~\ref{defn.pointwise.perm.action})}
\\
=&
\{(\pi\act x){\downarrow} \mid x\in\univi{n},\ \model{\phi[a\ssm x{\downarrow}{\uparrow}]}=\top\}
&\text{Lemma~\ref{lemm.arrows.equivariant}}
\\
=&
\{x{\downarrow} \mid x\in\univi{n},\ \model{\phi[a\ssm (\pi^\mone\act x){\downarrow}{\uparrow}]}=\top\}
&\text{Permutation action bij. on $\univi{n}$}
\\
=&
\{x{\downarrow} \mid x\in\univi{n},\ \model{\phi[a\ssm \pi^\mone\act (x{\downarrow}{\uparrow})]}=\top\}
&\text{Lemma~\ref{lemm.arrows.equivariant}}
\\
=&
\{x{\downarrow} \mid x\in\univi{n},\ \model{\phi[a\ssm \pi^\mone\mact (x{\downarrow}{\uparrow})]}=\top\}
&\text{Definition~\ref{defn.permutation.action}(1)}
\\
=&
\{x{\downarrow} \mid x\in\univi{n},\ \model{\phi[a\ssm x{\downarrow}{\uparrow}]}=\top\}
&\text{Corollary~\ref{corr.sub.pi.mact.x}},\ \pi^\mone\in\fix(\adjsupp(\phi))
\\
=& X
&\text{Fact}
\end{array}
$$
\end{proof}

\begin{nttn}
\label{nttn.extent}
Write $\f{extent}(\phi)$, read the \deffont{extent} of $\phi$, for the difference between the highest level variable or constant mentioned in $\phi$, and the lowest, plus $1$.

Thus for example if $\level(b)=10$ and $\level(a)=11$ then $\f{extent}(\tall a.\texi b.b\tin a)=2$.
By convention, if $\phi$ mentions no variables or constants then we set $\f{extent}(\phi)=0$. 
\end{nttn}

\begin{prop}[\bf Comprehension]
\label{prop.comprehension}
Suppose $\phi$ is a predicate with free variables in $\bar c$ and $a$ (so $\phi$ may not mention constants) and $\f{extent}(\phi)<n$ (where $n$ was chosen fixed but arbitrary in Definition~\ref{defn.n}).
Then 
$$
\model{\tall\bar c.\texi b.\tall a.a\tin b\tiff \phi}=\top.
$$
\end{prop}
\begin{proof}
Write $\level(a)=n\plus i$ for some $i\in\mathbb Z$.
Because $\f{extent}(\phi)<n$, 
we may use Proposition~\ref{prop.TA} to shift $\phi$ and suppose without loss of generality that $i=0$; thus $\level(a)=n$. 

For each variable symbol $c$ in $\bar c$ we choose a value $z\in\univi{\level(c)}$ and replace $c$ with $z$ in $\phi$ to obtain a predicate which we can write as $\phi[\bar c\ssm \bar z]$, which has free variables in $\{a\}$.

Write 
$$
C=\{x\in\univi{n} \mid \model{\phi[\bar c\ssm\bar z][a\ssm x]}=\top \}\subseteq\univi{n}.
$$
To check that $C\in\univi{n\plus 1}$ we must find an $X\in\univi{n}$ such that $C=X{\uparrow}\in\univi{n\plus 1}$.
We set 
$$
X=\{x{\downarrow} \mid x\in\univi{n},\ \model{\phi[\bar c\ssm\bar z][a\ssm x{\downarrow}{\uparrow}]}=\top \} \subseteq\univi{n\minus 1} .
$$
By Corollary~\ref{corr.X.in.Un} $X\in\univi{n}$.
Also, by Lemma~\ref{lemm.updown} $x\in C$ if and only if $x{\downarrow}{\uparrow}\in C$ and by construction this is if and only if $x{\downarrow}\in X$.
By Lemma~\ref{lemm.arrows.adjoint} this is if and only if $x\in X{\uparrow}$.
Thus $C=X{\uparrow}$.

The choice of $\bar z$ was arbitrary, so by Definition~\ref{defn.denotation} $\model{\tall\bar c.\texi b.\tall a.a\tin b\tiff \phi}=\top$ as required.
\end{proof}

\section{Consistency of TST+, and thus of NF}

\subsection{The main theorem}

\begin{thrm}
TST+ is consistent and thus --- since by a result of Specker~\cite{specker:typa} TST+ is equiconsistent with NF --- so is NF.
\end{thrm}
\begin{proof}
Suppose we have a derivation $\Delta$ of $\tbot$ in TST+ as defined in Definition~\ref{defn.nf.entailment}(2).
Let $n$ be the difference between the highest and the lowest degree variable in $\Delta$ (as per Notation~\ref{nttn.extent}), plus one.
Then $\univi{}$ clearly models first-order logic as per Definition~\ref{defn.denotation}, and by Propositions~\ref{prop.extensionality}, \ref{prop.comprehension}, and~\ref{prop.TA} we have that extensionality, comprehension (for all the instances that can appear in $\Delta$), and typical ambiguity are valid in it, contradicting our assumption that we have a derivation of $\tbot$.

We sketch the equiconsistency argument in Appendix~\ref{sect.tst.nf} for reference.
\end{proof}

\begin{rmrk}
\label{rmrk.consistency.strength}
The consistency strength of this proof is $\beth_\omega=\bigvee_{i\in\mathbb N}\beth_i=\bigvee_{i\in\mathbb N}\card\powerset^i(\mathbb N)$ (cf. Definition~\ref{defn.set.of.atoms} and Remark~\ref{rmrk.words.on.foundations}).
\end{rmrk}

\subsection{Elaborating on the model}
\label{subsect.elaborating.the.model}

Our constructions above deliver universes $\univi{}^n$ that validate comprehension of extent (Notation~\ref{nttn.extent}) up to some fixed but arbitrary $n$.
This suffices to prove consistency of TST+.

However, we can use some model theory to package up our constructions and deliver families of models of TST+/NF each member of which validates \emph{all} comprehensions, with no bound on extent.

We very briefly sketch two methods: one uses nonstandard models, and the other uses ultraproducts and {\L}o\'s' theorem.
What follows is just meant to note what we can do; the constructions are standard and a reader familiar with the model theory could fill in the rest.

\begin{rmrk}[Nonstandard models]
Let $\mathcal M$ be a nonstandard model of Zermelo set theory, by which we mean a model of Zermelo set theory whose natural numbers include an element which, when viewed from outside the model, is nonstandard~\cite[Figure~1]{kaye:modpa}.
Write this element $n$. 

We perform the constructions in this paper within that model, and obtain $\univi{n}^n$ for this (nonstandard) $n$.
Then any $\phi$ in TST+ can be easily interpreted in this model (e.g. a universal quantification $\tall a.\phi$ is interpreted to range over elements of the model at the suitable nonstandard level).
By easy facts of nonstandard models of arithmetic (clearly presented in Subsection~6.2 of~\cite{kaye:modpa}, especially Figure~3) we have infinite descending and ascending chains of nonstandard numbers $n\minus 1$, $n\plus 1$, $n\minus 2$, $n\plus 2$, \dots which we can use to interpret predicates of arbitrary extent.
There is no contradiction with our external notion of finiteness, because this is different from the model of finiteness inside $\mathcal M$. 
We can therefore write:
\end{rmrk}

\begin{thrm}
Any choice of nonstandard model $\mathcal M$ of Zermelo set theory, and nonstandard number $n\in\mathbb N_{\mathcal M}$ in that model, naturally yields a model of TST+, as $\univi{\mathcal M}^n$.
\end{thrm}

For the ultraproduct method, we can proceed as follows:
\begin{thrm}
Fix a nonprincipal ultrafilter $\mathcal U$ on $\mathbb N$.
We can give TST+ a semantics in which the denotation of level $i$ is the \emph{ultraproduct} $\hfs{i}/\mathcal U$ ($\mathbb N$-tuples of elements from $\hfs{i}$, quotiented such that two tuples are identified when they agree on a set of elements in $\mathcal U$).

For $y\in\hfs{i}/\mathcal U$ and $x\in\hfs{i\plus 1}/\mathcal U$, we set $\model{y\tin x}=\top$ when $\{n\in\mathbb N \mid \model{y_n\tin x_n}^n=\top\}\in\mathcal U$ (where $y_n\in\hfs{i}$ is the $n$th element of a representative tuple for $y$ --- recall that $y$ is an equivalence class of tuples --- and similarly for $x_n$).
The interpretation of $\tbot$, $\timp$, and $\tall$ over $\hfs{}/\mathcal U$ follows as standard for first-order syntax. 

Note that for any $\phi$, comprehension holds at sufficiently large $n$ --- because the extent of $\phi$ is fixed, so as we increase $n$, $n$ will eventually exceed the extent of $\phi$.
It follows that this model satisfies comprehension.
The rest follows by standard methods: the interested reader can look up {\L}o\'s' theorem \cite[page 90]{bell:modu}. 
\end{thrm}

\hyphenation{Mathe-ma-ti-sche}
\providecommand{\bysame}{\leavevmode\hbox to3em{\hrulefill}\thinspace}
\providecommand{\MR}{\relax\ifhmode\unskip\space\fi MR }
\providecommand{\MRhref}[2]{%
  \href{http://www.ams.org/mathscinet-getitem?mr=#1}{#2}
}
\providecommand{\href}[2]{#2}

\appendix

\section{Appendix}

\subsection{Sketch of transformation between (countable) models of TST+ and NF}
\label{sect.tst.nf}

In this Appendix we will re-express in a slightly more modern notation the \emph{simple model-theoretic argument} described towards the end of the third page of~\cite{specker:typa}.
First, a preliminary definition:
\begin{defn}
A \deffont{countably infinite model of Quine's NF} consists of a countably infinite set $V$ and a \deffont{membership relation} $\nfin \subseteq V\times V$, such that extensionality and (stratifiable) comprehension are valid over the structure $(V,\nfin)$.
\end{defn}

\begin{thrm}
\label{thrm.specker}
There exists a countably infinite model of TST+ if and only if there exists a countably infinite model of Quine's NF.
\end{thrm}
\begin{proof}
We outline the manipulation on models from~\cite{specker:typa} to convert between the models:

Suppose we have a countably infinite sets model of TST+, so we have 
$$
(V_0, V_1\subseteq\powerset(V_0), V_2\subseteq\powerset(V_1), \dots)
$$
and each $V_i$ is countably infinite.

To build a model of Quine's NF, we choose isomorphisms $f_i:V_i\cong V_{i\plus 1}$ and set 
$$
V=V_0
\quad\text{and}\quad
y\nfin x
\ \text{when}\ 
y\in f_i(x) .
$$
In this way, $f_i$ assigns a $V_i$-extension to each element of $V_i$.
It is now routine to check that $(V,\nfin)$ is an extensional model (because the sets model is) with stratifiable comprehension (because the sets model has stratified comprehension) and so is a model of Quine's NF.

Conversely, given a countable model $(V,\nfin)$ of Quine's NF we construct a countable model of TST+ by setting $V_0=V$ and $V_1=\{ \{v'{\in}V \mid v'\nfin v \} \mid v\in V\}$ and repeatedly unrolling the sets extension to obtain a sequence of countably infinite sets $V_i\subseteq\powerset^i(V)$.
\end{proof}

\subsection{Generalising Russell's paradox: cyclic stratified syntax does not normalise}
\label{sect.aside}

The comprehension rewrite rule on TST syntax augmented with stratified comprehension terms $\stxst{a}{\phi}$ --- a.k.a. $\beta$-reduction 
$$
t\tin\stxst{a}{\phi}\longrightarrow \phi[a\ssm t]  
$$
--- is normalising~\cite{gabbay:lanssc}. 
What if we stratify, but using some cyclic graph like integers modulo some $n$, instead of $\mathbb N$?
We take a few moments to illustrate that we do indeed lose normalisation:

\begin{defn}
Define a na\"ive sets term $r_1$ by
$$
r_1=\stxst{a}{\tneg (a\tin a)}.
$$
\end{defn}

\begin{rmrk}
\emph{Russell's paradox} corresponds to the following rewrite:
$$
r_1\tin r_1 \to \tneg (r_1\tin r_1) .
$$
This rewrite is non-normalising, in the sense that the term we started with is a subterm of the term we finish with --- we have a loop.
\end{rmrk}

\begin{nttn}
Given a na\"ive sets term $t$ and $i\geq 0$, define $\iota^i t$ by:
$$
\iota^0 t = t
\quad\text{and}\quad
\iota t = \stxst{a}{t\tin a}
\quad\text{and}\quad
\iota^{i\plus 1} t = \iota(\iota^i t) 
.
$$
\end{nttn}

\begin{rmrk}
We note the following rewrites: 
$$
s\tin \iota t \to t\tin s
\quad\text{and}\quad
\iota t\tin\iota s \to s\tin\iota t \to t\tin s .
$$
\end{rmrk}

\begin{defn}
We define
$$
r_2 = \stxst{a}{\tneg(\iota a\tin a)}.
$$
\end{defn}

\begin{rmrk}
If we try to normalise $\iota r_2\tin r_2$, we obtain:
$$
\iota r_2\tin r_2 
\to
\tneg(\iota^2 r_2\tin\iota r_2)
\to^*
\tneg(\iota r_2\tin r_2) .
$$
($\to^*$ above denotes multiple rewrites.) 
\end{rmrk}

\begin{defn}
Generalising the above, if $n\geq 1$ then we write
$$
r_n = \stxst{a}{\tneg(\iota^{n\minus 1} a\tin a)}.
$$
\end{defn}

\begin{rmrk}
If we try to normalise $\iota^{n\minus 1}r_n\tin r_n$, we obtain
$$
\iota^{n\minus 1}r_n\tin r_n 
\to
\tneg(\iota^{2n\minus 2}r_n\tin\iota^{n\minus 1}r_n)
\to^*
\tneg(\iota^{n\minus 1}r_n\tin r_n) 
.
$$
\end{rmrk}

So if we stratify na\"ive sets with finitely many levels $\mathbb Z_{\text{mod}(n)}$ for any finite $n$, then there exist predicates with no normal form.

\end{document}